\documentclass[10pt]{amsart}
\usepackage{latexsym}
\usepackage{amssymb}
\usepackage{amsmath}
\usepackage{mathrsfs}
\usepackage{bbm}
\usepackage{bbold}

\def\ssection#1{\bigskip\begin{center}{\scshape#1}\end{center}\medskip}
\def\subsection#1{\smallskip\par{\bf #1.}}
\def\subsubsec#1{{\bf#1}}

 \def\roman#1{\mathrm{#1}}
 \def\olim{\mathop{o\text{-}\!\lim}}
 \def\bolim{\mathop{bo\text{-}\!\lim}}
 \def\dom{\mathop{\fam0 dom}}
 \def\Orth{\mathop{\fam0 Orth}}
 \def\On{\mathop{\fam0 On}}
 \def\im{\mathop{\fam0 im}}
 \def\ZFC{\mathop{\fam0 ZFC}}
 \def\End{\mathop{\fam0 End}}
 \def\mix{\mathop{\fam0 mix}}
 \def\sa{\mathop{\fam0 sa}}
 \def\Id{\mathop{\fam0 Id}}
 \def\Fin{\mathop{\fam0 fin}\nolimits}
 \def\osum{o\text{-}\!\sum}
 \def\bosum{bo\text{-}\!\sum}
 \def\Re{\mathop{\fam0 Re}}

\def\upwardarrow{\mathord{
  \hbox to 5pt{\hss$\vcenter{\hbox to 2.4pt{\hss$\mathchar"222$\hss}\hrule}\hss$}
}}
\def\downwardarrow{\mathord{
  \hbox to 5pt{\hss$\vcenter{\hrule\hbox to 2.4pt{\hss$\mathchar"223$\hss}}\hss$}
}}

 \def\endproc{\rm}
 \def\proclaim#1{{\bf #1}\sl}
 \def\Proclaim#1{\smallskip\par{\bf #1}\sl}
\def\tvert{|\!|\!|}
\def\[{\mathopen{\kern1pt \vrule height7.5pt depth1.5pt width1pt\kern1.5pt}}
\def\]{\mathclose{\kern1pt \vrule height7.5pt depth1.5pt width1pt\kern1pt}}
\begin{document}
\title[Boolean Methods]{Boolean Methods\\ in the Theory of Vector Lattices}
\begin{abstract}
This is an overview of the recent results of interaction
of Boolean valued analysis and vector lattice theory.
\end{abstract}
\address[]{\noindent
\vbox{
Institute of Applied Mathematics and Informatics\hfill\break
\indent Vladikavkaz, RUSSIA}
\vskip 5truept
\indent\vbox{
Sobolev Institute of Mathematics\hfill\break
\indent Novosibirsk, RUSSIA
}}
\email{
kusraev@alanianet.ru\\
sskut@math.nsc.ru
}
\author{A. G. Kusraev and S. S. Kutateladze}
%\date{July 26, 2005}

\maketitle

%\baselineskip=.95\baselineskip
\ssection{Introduction}

Boolean valued analysis is a~general mathematical method that rests on
a~special model-theoretic technique. This technique consists  primarily
in comparison between the representations of  arbitrary mathematical objects
and theorems in two different set-theoretic models whose constructions start
with principally distinct Boolean algebras.  We usually take as these models
the cosiest Cantorian paradise, the von Neumann universe  of
Zermelo--Fraenkel set theory, and a~special universe of Boolean valued
``variable'' sets trimmed and chosen so that the traditional concepts and
facts of mathematics acquire completely unexpected and bizarre
interpretations.  The use of two models, one of which is {\it formally\/}
nonstandard, is a~family feature of {\it nonstandard  analysis}.  For this
reason, Boolean valued analysis  means an instance of nonstandard analysis in
common parlance.  By the way, the term {\it Boolean valued analysis\/} was
minted by G.~Takeuti.

Proliferation of Boolean valued models is due to P.~Cohen's final
breakthrough in Hilbert's Problem Number One.  His method of forcing was
rather intricate and the inevitable attempts at simplification gave rise to
the Boolean valued models by  D.~Scott,  R.~Solovay, and P.~Vop\v enka.

Our starting point is a brief description of the best Cantorian paradise in
shape of the von Neumann universe and a specially-trimmed Boolean valued
universe that are usually taken as these two models. Then we present a
special ascending and descending machinery for interplay between the models.
We consider the reals and complexes inside a Boolean valued model by using
the celebrated Gordon's Theorem which we read as follows: Every universally
complete vector lattice is an interpretation of the reals in an appropriate
Boolean-valued model.  We proceed with demonstrating the Boolean valued
approach to the two familiar problems:  (1) When is a band preserving operator
order bounded?  (2) When is an order bounded operator a sum or difference of
two lattice homomorphisms? In conclusion we briefly overview
the details of some typical spaces and operators
together with their Boolean valued representations.

\ssection {1. Boolean  Requisites}

We start with recalling some auxiliary facts about the construction
and treatment of the von Neumann universe and a specially-trimmed
Boolean valued universe.

\subsection{1.1}~The {\it von Neumann universe} ${\mathbbm V}$ results
by transfinite recursion over ordinals.
As the initial object of this construction we take the
empty set.  The elementary step of introducing new sets consists
in  uniting the powersets of the sets already
available. Transfinitely repeating these steps, we exhaust
the class of all sets. More precisely, we put
$\mathbbm V\!:=\cup_{\alpha\in \operatorname{On}}{\mathbbm V}_{\alpha}$,
where $\operatorname{On}$ is the class of all ordinals and
$$
\gathered
{\mathbbm V}_0\!:=\varnothing,
\\
{\mathbbm V}_{\alpha+1}\!:={\mathscr P}({\mathbbm V}_{\alpha}),
\\
{\mathbbm V}_{\beta}\!:=\bigcup\limits_{\alpha<\beta}{\mathbbm V}_{\alpha}
\quad (\text{$\beta$ is a limit ordinal}).
\endgathered
$$
The class~${\mathbbm V}$ is the standard model of
Zermelo--Fraenkel set theory.

\subsection{1.2}
% Universes and truth values}
Let ${\mathbbm B}$ be a~complete Boolean algebra.
Given an ordinal $\alpha$, put
$$
%\begin{displaymath}
\gathered
{\mathbbm V}_{\alpha}^{({\mathbbm B})}\!:=
\{x: x \mbox{ is\ a\ function}\ \\
\wedge\ (\exists
\beta)(\beta<\alpha\ \wedge \dom (x)
\subset {\mathbbm V}_{\beta}^{({\mathbbm B})}\ \wedge\
\im (x)\subset {\mathbbm B})\}.
\endgathered
%\end{displaymath}
$$

After this recursive definition the {\it Boolean valued
universe\/}
${\mathbbm V}^{({\mathbbm B})}$
or, in other words, the {\it class of ${\mathbbm B}$-sets\/}
is introduced by
$$
{\mathbbm V}^{({\mathbbm B})}\!:=\bigcup\limits_{\alpha\in\On} {\mathbbm V}_{\alpha}^{({\mathbbm B})},
$$
with $\On$ standing for the class of all ordinals.

In case of the two-element Boolean algebra
$\mathbb 2:=\{\mathbb 0, \mathbb 1\}$
this procedure yields a~version
of the classical
von Neumann universe ${\mathbbm V}$.

Let $\varphi$ be an arbitrary formula of $\ZFC$,
Zermelo--Fraenkel set theory with  choice.
The {\it Boolean truth value\/}
$[\![\varphi]\!]\in {\mathbbm B}$
is introduced by induction on the length of a~formula~$\varphi$
by using the natural interpretation of  the propositional connectives and
quantifiers in~${\mathbbm B}$
and  the way in which $\varphi$ results from atomic formulas. The Boolean truth values of the
{\it atomic formulas\/}
$x\in y$
and
$x=y$,
with
$x,y\in {\mathbbm V}^{({\mathbbm B})}$,
are defined by means of the following recursion schema:
$$
\gathered
{}
[\![x\in y]\!]=
\bigvee\limits_{t\in\dom(y)} y(t)\wedge [\![t=x]\!],
\\[0.5\jot]
[\![x=y]\!]=\bigvee\limits_{t\in\dom(x)} x(t)\Rightarrow [\![t\in
y]\!]\wedge \bigvee\limits_{t\in\dom(y)} y(t)\Rightarrow [\![t\in
x]\!].
\endgathered
$$
The sign
$\Rightarrow$
symbolizes the implication in
${\mathbbm B}$; i.e.,
$a\Rightarrow b:=a^\ast\vee b$
where
$a^\ast$
is as usual the {\it complement\/} of~$a$.

%\subsection{Principles of analysis}
The universe ${\mathbbm V}^{({\mathbbm B})}$
with the  Boolean truth value of a~formula is a~model
of set theory in the sense that every theorem of $\ZFC$
is true inside ${\mathbbm V}^{({\mathbbm B})}$.

\subsection{1.3}
\proclaim{Transfer Principle.}
For every theorem
$\varphi$ of~$\ZFC$, we have
 $[\![\varphi]\!]=\mathbb 1$;
i.e.,
$\varphi$
is true inside
${\mathbbm V}^{({\mathbbm B})}$.
\endproc

Enter into the next agreement: If
$x$
is an element of
${\mathbbm V}^{({\mathbbm B})}$
and
$\varphi (\cdot)$
is a~formula of $\ZFC$, then the phrase
``$x$
satisfies
$\varphi$
inside ${\mathbbm V}^{({\mathbbm B})}$''
or, briefly,
``$\varphi (x)$
is true inside
${\mathbbm V}^{({\mathbbm B})}$'' means that
$[\![\varphi(x)]\!]=\mathbb  1$. This is sometimes written
as ${\mathbbm V}^{({\mathbbm B})}\models \varphi (x)$.

Given $x\in {\mathbbm V}^{({\mathbbm B})}$ and $b\in {\mathbbm B}$,
define the function $bx:z\mapsto bx (z)\quad (z\in\dom (x))$. Here
we presume that $b\varnothing\!:=\varnothing$ for all $b\in {\mathbbm B}$.

There is a~natural equivalence relation $x\sim y\leftrightarrow
[\![x=y]\!]=\mathbb 1$ in  the class ${\mathbbm V}^{({\mathbbm B})}$.
Choosing a~representative of the smallest rank in each equivalence class or,
more exactly, using the so-called ``Frege--Russell--Scott trick,'' we obtain
a~{\it separated Boolean valued universe\/} $\overline {\mathbbm V}^{{\,}
({\mathbbm B})}$ in which $ x=y\leftrightarrow [\![x=y]\!]=\mathbb 1$.

It is easily to see that the Boolean truth value of a~formula remains
unaltered if we replace in it each element of ${\mathbbm V}^{({\mathbbm B})}$
by one of its equivalents. In this connection from now on we take ${\mathbbm
V}^{({\mathbbm B})}:=\overline {\mathbbm V}^{{\,} ({\mathbbm B})}$ without
further specification.

Observe that in $\overline {\mathbbm V}^{{\,} ({\mathbbm B})}$ the element
$bx$ is defined correctly for $x\in\overline {\mathbbm V}^{{\,} ({\mathbbm
B})}$ and $b\in {\mathbbm B}$ since $[\![x_1=x_2]\!]={\mathbb  1}\rightarrow
[\![bx_1=b x_2]\!]= b\Rightarrow [\![x_1=x_2]\!]={\mathbb  1}$.  For
a~similar reason, we often write ${\mathbb  0}:=\varnothing$, and in
particular ${\mathbb  0}\varnothing =\varnothing={\mathbb  0} x$ for $x\in
{\mathbbm V}^{({\mathbbm B})}$.

\subsection{1.4}
\proclaim{Mixing Principle.}
Let
$(b_{\xi})_{\xi\in\Xi}$
be a~{\it partition of unity\/} in~
${\mathbbm B}$,
i.e.
$\sup_{\xi\in\Xi} b_{\xi}$$ =\sup {\mathbbm B}=\mathbb  1$
and
$\xi\neq\eta\rightarrow b_{\xi}\wedge b_{\eta}=\mathbb 0$.
To each family
$(x_{\xi})_{\xi\in\Xi}$
in~${\mathbbm V}^{({\mathbbm B})}$
there exists a~unique element
$x$
in the separated universe such that
$[\![x=x_{\xi}]\!]\ge b_{\xi}\quad (\xi\in\Xi).
$
\endproc

This element is called the {\it mixing\/} of $(x_{\xi})_{\xi\in\Xi}$
by~$(b_{\xi})_{\xi\in\Xi}$ and is denoted by $\sum_{\xi\in\Xi} b_{\xi}
x_{\xi}$. Thus, the mixing principle asserts that every Boolean valued
universe is rich in mixings.

\subsection{1.5}
\proclaim{Maximum Principle.}
The least upper bound is attained on the right-hand side of the
formula for the Boolean truth-value of the existential quantifier.
More precisely, if $\varphi$ is a~formula of~$\ZFC$ then there
is a ${\mathbbm B}$-valued set $x_0$ satisfying
$[\![(\exists x)\varphi (x)]\!]=[\![\varphi (x_0)]\!]$.
\endproc

\ssection{2. The Escher Rules}

Boolean valued analysis  consists primarily in comparison of the
instances of a~mathematical object or idea in two  Boolean
valued models. This is impossible to achieve without some
dialog  between the universes~${\mathbbm V}$ and~${\mathbbm
V}^{({\mathbbm B})}$.  In other words, we need a~smooth
mathematical toolkit for revealing interplay between the
interpretations of one and the same fact in the two
models~${\mathbbm V}$ and~${\mathbbm V}^{({\mathbbm B})}$.  The
relevant {\it ascending-and-descending technique\/} rests on the
functors of canonical embedding, descent, and ascent.

\subsection{2.1}
%Embedding}
We start with the canonical embedding of the von
Neumann universe $\mathbbm V$.

Given
${x\in {\mathbbm V}}$,
we denote by~$x^{\scriptscriptstyle\wedge}$
the {\it standard name\/}
of~$x$
in~${\mathbbm V}^{({\mathbbm B})}$;
i.e., the element defined by the following recursion schema:
$
\varnothing^{\scriptscriptstyle\wedge}\!:=\varnothing$,\quad
$\dom (x^{\scriptscriptstyle\wedge}):=
\{y^{\scriptscriptstyle\wedge} : y\in x \}$,\quad
$\im (x^{\scriptscriptstyle\wedge})\!:=\{{\mathbb 1}\}$.
Observe some properties of the mapping
$x\mapsto x^{\scriptscriptstyle\wedge}$
we need in the sequel.

{\bf(1)}
For an~arbitrary
$x\in {\mathbbm V}$
and a~formula
$\varphi$
of~$\ZFC$ we have
$$
\gathered
{}[\![(\exists y\in
x^{\scriptscriptstyle\wedge})\,\varphi(y)]\!]=
\bigvee\limits_{z\in x}
[\![\varphi(z^{\scriptscriptstyle\wedge})]\!],
\\
[\![(\forall y\in x^{\scriptscriptstyle\wedge})\,\varphi(y)]\!]
=\bigwedge\limits_{z\in x}[\![\varphi(z^{\scriptscriptstyle\wedge})]\!].
\endgathered
$$

{\bf(2)}
If
$x$ and~$y$
are elements of~${\mathbbm V}$
then, by transfinite induction, we establish
$x\in y\leftrightarrow {\mathbbm V}^{({\mathbbm B})}\models x^{\scriptscriptstyle\wedge}\in y^{\scriptscriptstyle\wedge},
\quad
x= y\leftrightarrow {\mathbbm V}^{({\mathbbm B})}\models x^{\scriptscriptstyle\wedge}=y^{\scriptscriptstyle\wedge}$.
In other words, the standard name can be considered
as an~embedding of~${\mathbbm V}$
into~${\mathbbm V}^{({\mathbbm B})}$.
Moreover, it is beyond a~doubt that the standard name sends
${\mathbbm V}$
onto
${\mathbbm V}^{(\mathbb 2)}$,
which fact is demonstrated by the next proposition:

{\bf(3)}
The following holds:
$
(\forall u\in {\mathbbm V}^{(\mathbb 2)})\,(\exists ! x\in {\mathbbm V})\
{\mathbbm V}^{({\mathbbm B})}\models u=x^{\scriptscriptstyle\wedge}.
$

A~formula is called {\it bounded\/}
or {\it restricted\/}
if each bound variable in it is restricted by a~bounded
quantifier; i.e., a~quantifier ranging over a~particular set.
The latter means that each bound variable~$x$
is restricted by a~quantifier of the form~$(\forall x\in y)$
or
$(\exists x\in y)$
for some~$y$.

\subsection{2.2}
\proclaim{Restricted Transfer Principle.}
For each bounded formula~$\varphi$
of~$\ZFC$ and every collection~$x_1,\dots,x_n\in {\mathbbm V}$
the following  holds:
$
\varphi(x_1,\dots,x_n)\leftrightarrow {\mathbbm V}^{({\mathbbm B})}\models
\varphi(x_1^{\scriptscriptstyle\wedge},\dots,x_n^{\scriptscriptstyle\wedge}).
$
\endproc
Henceforth, working in the separated universe~$\overline{{\mathbbm V}}^{({\mathbbm B})}$,
we agree to preserve the symbol~$x^{\scriptscriptstyle\wedge}$
for the distinguished element of the class
corresponding to~$x$.

Observe for example that
 the restricted transfer principle yields:
$$
\gathered
\text{``$\Phi$ is a correspondence from $x$ to
$y$''}
%\leftrightarrow
\\
\leftrightarrow
{\mathbbm V}^{({\mathbbm B})}\models
\text{``$\Phi^{\scriptscriptstyle\wedge}$
is a correspondence from $x^{\scriptscriptstyle\wedge}$ to
$y^{\scriptscriptstyle\wedge}$'';}
\\
\text{``$f: x\to y$''
$\leftrightarrow $
${\mathbbm V}^{({\mathbbm B})}\models$
``$f^{\scriptscriptstyle\wedge}: x^{\scriptscriptstyle\wedge}\to
y^{\scriptscriptstyle\wedge}$''}
\endgathered
$$
\noindent
(moreover,
$f(a)^{\scriptscriptstyle\wedge}=f^{\scriptscriptstyle\wedge}(a^{\scriptscriptstyle\wedge})$
for all $a\in x$).
Thus, the standard name can be considered as a
covariant functor of the category of sets (or
correspondences) inside~${\mathbbm V}$
to an~appropriate subcategory of~${\mathbbm V}^{(\mathbb 2)}$
in the separated universe~${\mathbbm V}^{({\mathbbm B})}$.

\subsection{2.3}~A~set
$X$
is  {\it finite\/}
if
$X$
coincides with the image of a~function on a~finite ordinal.
In symbols, this is expressed as
$\Fin (X)$; hence,
$$
\Fin(X):=(\exists\,n)(\exists\,f)(n\in\omega\wedge
f \mbox{ is a function} \wedge\ \dom(f)=n\wedge\im(f)=X)
$$
(as usual $\omega:=\{0,1,2,\dots\}$).
Obviously, the above formula is  not bounded.
Nevertheless there is a~simple transformation rule for the  class
of finite sets under the canonical embedding.
Denote by
${\mathscr P}_{\Fin}(X)$
the class of all finite subsets of $X$; i.e.,
${\mathscr P}_{\Fin}(X):=\{Y\in{\mathscr P}(X):\Fin (Y)\}.
$
For an arbitrary set $X$ the following holds:
$
{\mathbbm V} ^{({\mathbbm B})}\models\mathscr P _{\Fin} (X)^{\scriptscriptstyle\wedge} =
\mathscr P _{\Fin} (X^{\scriptscriptstyle\wedge}).
$

\subsection{2.4}
%. Descent}
Given an~arbitrary element~$x$
of the (separated) Boolean  valued universe~${\mathbbm V}^{({\mathbbm B})}$,
we define the {\it descent}
$x{\downarrow}$
of~$x$ as $
x{\downarrow}:=\{y\in {\mathbbm V}^{({\mathbbm B})}: [\![y\in x]\!]={\mathbb 1} \}$.
We list the simplest properties of descending:

{\bf(1)}
The class~$x{\downarrow}$
is a~set, i.e.,
$x{\downarrow}\in {\mathbbm V}$
for all~$x\in {\mathbbm V}^{({\mathbbm B})}$.
If
$[\![x\ne \varnothing]\!]={\mathbb 1}$
then
$x{\downarrow}$
is a~nonempty set.

{\bf(2)}
Let
$z\in {\mathbbm V}^{({\mathbbm B})}$
and
$[\![z\ne \varnothing]\!]={\mathbb 1}$.
Then for every formula~$\varphi$
of~$\ZFC$ we have
$$
\gathered
{}
[\![(\forall x\in z)\,\varphi(x)]\!]
=\bigwedge\limits_{ x\in z{\downarrow}}[\![\varphi(x)]\!],
\\
[\![(\exists x\in z)\,\varphi(x)]\!]
=\bigvee\limits_{x\in z{\downarrow}}[\![\varphi(x)]\!].
\endgathered
$$
Moreover, there exists
$x_0\in z{\downarrow}$
such that
$[\![\varphi(x_0)]\!]=[\![(\exists x\in z)\,\varphi(x)]\!]$.

{\bf(3)}
Let
$\Phi$
be a~correspondence from~$X$
to~$Y$
in
${\mathbbm V}^{({\mathbbm B})}$.
Thus,
$\Phi$,
$X$,
and~$Y$
are elements of~${\mathbbm V}^{({\mathbbm B})}$
and, moreover,
$[\![\Phi\subset X\times Y]\!]={\mathbb 1}$.
There is a~unique correspondence~$\Phi{\downarrow}$
from
$X{\downarrow}$
to~$Y{\downarrow}$
such that
$
\Phi{\downarrow}(A{\downarrow})=\Phi(A){\downarrow}
$
for every nonempty subset~$A$
of~$X$
inside~${\mathbbm V}^{({\mathbbm B})}$.
The correspondence~$\Phi{\downarrow}$
from
$X{\downarrow}$
to
$Y{\downarrow}$
of the above proposition is called the
{\it descent\/}
of the correspondence~$\Phi$
from~$X$
to~$Y$
inside~${\mathbbm V}^{({\mathbbm B})}$.

{\bf(4)}
The descent of the composite of correspondences inside~${\mathbbm V}^{({\mathbbm B})}$
is the composite of their descents: $(\Psi\circ\Phi){\downarrow}=\Psi{\downarrow}\circ\Phi{\downarrow}.
$

{\bf(5)}
If
$\Phi$
is a~correspondence inside~${\mathbbm V}^{({\mathbbm B})}$
then
$
(\Phi^{-1}){\downarrow}=(\Phi{\downarrow})^{-1}.
$

{\bf(6)}
Let
$\Id_X$
be the identity mapping inside~${\mathbbm V}^{({\mathbbm B})}$
of a~set~$X\in {\mathbbm V}^{({\mathbbm B})}$.
Then
$
({\Id}_X){\downarrow}={\Id}_{X{\downarrow}}.
$

{\bf(7)}
Suppose that
$X,Y,f\in {\mathbbm V}^{({\mathbbm B})}$
are such that
$[\![f:X\to Y]\!]={\mathbb 1}$,
i.e.,
$f$
is a~mapping from~$X$
to~$Y$
inside~${\mathbbm V}^{({\mathbbm B})}$.
Then
$f{\downarrow}$
is a~unique mapping from~$X{\downarrow}$
to~$Y{\downarrow}$
satisfying
$[\![f{\downarrow}(x)=f(x)]\!]={\mathbb 1}$ for all $x\in X{\downarrow}$.

By virtue of~\hbox{(1)--(7)}, we can
consider the descent operation as a~functor from the category
of~${\mathbbm B}$-valued sets and mappings (correspondences)
to the category of the usual
sets and mappings (correspondences) (i.e., in the sense of~${\mathbbm V}$).

{\bf(8)}
Given
$x_1,\dots,x_n\in {\mathbbm V}^{({\mathbbm B})}$,
denote by
$(x_1,\dots,x_n)^{\mathbbm B}$
the corresponding ordered $n$-tuple inside~${\mathbbm V}^{({\mathbbm B})}$.
Assume that
$P$
is an~$n$-ary relation on~$X$
inside~${\mathbbm V}^{({\mathbbm B})}$;
i.e.,
$X,P\in {\mathbbm V}^{({\mathbbm B})}$
and
$[\![P\subset X^{n^{\scriptscriptstyle\wedge}}]\!]={\mathbb 1}$,
where $n\in \omega$.
Then there exists an~$n$-ary relation~$P'$
on~$X{\downarrow}$
such that
$(x_1,\dots,x_n)\in P'\leftrightarrow [\![(x_1,\dots,x_n)^{\mathbbm B}\in P]\!]={\mathbb 1}$.
Slightly abusing notation, we denote the relation~$P'$
by the same symbol~$P{\downarrow}$
and call it the {\it descent\/} of~$P$.

\subsection{2.5}
%. Ascent}
Let
$x\in {\mathbbm V}$
and
$x\subset {\mathbbm V}^{({\mathbbm B})}$;
i.e., let
$x$
be some set composed of
${\mathbbm B}$-valued sets or, in other words,
$x\in \mathscr P({\mathbbm V}^{({\mathbbm B})})$.
Put
$\varnothing{\uparrow}:=\varnothing$
and
$
\dom (x{\uparrow}):=x,\quad \im (x{\uparrow}):=\{{\mathbb 1}\}
$
if
$x\neq \varnothing$.
The element
$x{\uparrow}$
(of the separated universe~${\mathbbm V}^{({\mathbbm B})}$,
i.e., the distinguished representative of the
class~$\{y\in {\mathbbm V}^{({\mathbbm B})}: [\![y=x{\uparrow}]\!]={\mathbb 1}
 \}$)
is called the {\it ascent\/}
of~$x$.

{\bf(1)}
For all
$x\in \mathscr P({\mathbbm V}^{({\mathbbm B})})$
and every formula~$\varphi$ we have the following:
$$
\gathered
{}[\![(\forall z\in x{\uparrow})\,\varphi(z)]\!]=
\bigwedge_{y\in x} [\![\varphi (y)]\!],
\\
[\![(\exists z\in x{\uparrow})\,\varphi(z)]\!]=
\bigvee_{y\in x} [\![\varphi (y)]\!].
\endgathered
$$

Introducing the ascent of a~correspondence~$\Phi\subset X\times Y$,
we have to bear in mind a~possible
distinction between the domain of departure~$X$
and the domain
$\dom (\Phi):=\{x\in X:\Phi(x)\ne \varnothing \}$.
This circumstance is immaterial for the sequel;
therefore,  speaking of ascents, we
always imply total correspondences;
i.e.,
$\dom (\Phi)=X$.

{\bf(2)}
Let
$X,Y,\Phi\in {\mathbbm V}^{({\mathbbm B})}$,
and let
$\Phi$
be a~correspondence from~$X$
to~$Y$.
There exists a~unique correspondence~$\Phi{\uparrow}$
from~$X{\uparrow}$
to~$Y{\uparrow}$
inside~${\mathbbm V}^{({\mathbbm B})}$
such that
$
\Phi{\uparrow}(A{\uparrow})=\Phi(A){\uparrow}
$
is valid for every subset~$A$
of~$\dom (\Phi)$
if and only if
$\Phi$
is {\it extensional};
i.e., satisfies the condition
$
y_1\in \Phi(x_1)\to [\![x_1=x_2]\!]\le
\bigvee\nolimits_{y_2\in \Phi(x_2)} [\![y_1=y_2]\!]
$
for
$x_1,x_2\in \dom (\Phi)$.
In this event,
$\Phi{\uparrow}=\Phi'{\uparrow}$,
where
$\Phi':=\{(x,y)^{\mathbbm B}: (x,y)\in \Phi \}$.
The element
$\Phi{\uparrow}$
is called the {\it ascent\/}
of~$\Phi$.

{\bf(3)}
The composite of extensional correspondences is
extensional. Moreover, the ascent of a~composite
is equal to the composite of the ascents
inside~${\mathbbm V}^{({\mathbbm B})}$:
On assuming that
$\dom (\Psi)\supset \im (\Phi)$
we have
$
{\mathbbm V}^{({\mathbbm B})}\vDash(\Psi\circ\Phi){\uparrow}=\Psi{\uparrow}\circ\Phi{\uparrow}.
$

Note that if
$\Phi$
and
$\Phi^{-1}$
are extensional then
$(\Phi{\uparrow})^{-1}=(\Phi^{-1}){\uparrow}$.
However, in general, the extensionality of~$\Phi$
in no way guarantees the extensionality of~$\Phi^{-1}$.

{\bf(4)}
It is worth mentioning that if an~extensional correspondence~$f$
is a~function from~$X$
to~$Y$
then the ascent~$f{\uparrow}$ of~$f$
is a~function from~$X{\uparrow}$
to~$Y{\uparrow}$.
Moreover, the extensionality property can be stated
as follows:
$[\![x_1=x_2]\!]\le [\![f(x_1)=f(x_2)]\!]$  for all $x_1,x_2\in X$.

\subsection{2.6}
Given a~set
$X\subset {\mathbbm V}^{({\mathbbm B})}$,
we denote by the symbol~$\mix(X)$
the set of all mixings of the form~$\mix (b_{\xi}x_{\xi})$,
where
$(x_{\xi})\subset X$
and
$(b_{\xi})$
is an~arbitrary partition of unity.
The following propositions are referred to as
the {\it arrow cancellation rules\/}
or {\it ascending-and-descending rules}.
There are many good reasons to call them simply
the {\it Escher rules\/}~\cite{Hof}.

{\bf(1)}
Let
$X$
and
$X'$
be subsets of~${\mathbbm V}^{({\mathbbm B})}$
and let
$f:X\to X'$
be an~extensional mapping. Suppose that
$Y,Y',g\in {\mathbbm V}^{({\mathbbm B})}$
are such that
$[\![\,Y\ne \varnothing]\!]=[\![\,g:Y\to Y']\!]={\mathbb 1}$.
Then
$
%\gathered
X{\uparrow}{\downarrow}=\mix(X),
Y{\downarrow}{\uparrow}=Y,
%\\
f{\uparrow}{\downarrow}=f,
$
and
$
g{\downarrow}{\uparrow}=g.
%\endgathered
$

{\bf(2)}~From 2.3\,(8) we easily infer the useful relation:
 $
{\mathscr P}_{\Fin}(X{\uparrow})=\{{\theta{\uparrow}} :
\theta\in{\mathscr P}_{\Fin}(X)\}{\uparrow}.
 $

%\ssection{The technique of ascending and descending}
Suppose that
$X\in {\mathbbm V}$,
$X\ne\varnothing$;
i.e.,
$X$
is a~nonempty set. Let the letter~$\iota$
denote the standard name embedding
$x\mapsto x^{\scriptscriptstyle\wedge}$
$(x\in X)$.
Then
$\iota(X){\uparrow}=X^{\scriptscriptstyle\wedge}$
and
$X=\iota^{-1}(X^{\scriptscriptstyle\wedge}{\downarrow})$.
Using the above relations, we may extend the
descent and ascent operations to the case in which
$\Phi$
is a~correspondence from~$X$
to~$Y{\downarrow}$
and
$[\![\Psi$ is a~correspondence from~$X^{\scriptscriptstyle\wedge}$
to~$Y]\!]={\mathbb 1}$,
where
$Y\in {\mathbbm V}^{({\mathbbm B})}$.
Namely, we put
$\Phi\upwardarrow:=(\Phi\circ\iota){\uparrow}$
and
$\Psi\downwardarrow:=\Psi{\downarrow}\circ\iota$.
In this case,
$\Phi\upwardarrow$
is called the
{\it modified ascent\/}
of~$\Phi$
and
$\Psi\downwardarrow$
is called the
{\it modified descent\/}
of~$\Psi$.
(If the context excludes ambiguity then we
briefly speak of ascents and descents using simple
arrows.)
It is easy to see that
$\Psi\upwardarrow$
is a~unique correspondence inside~${\mathbbm V}^{({\mathbbm B})}$
satisfying the relation
$
[\![\Phi\upwardarrow(x^{\scriptscriptstyle\wedge})=\Phi(x){\uparrow}]\!]={\mathbb 1}\quad (x\in X).
$
Similarly,
$\Psi\downwardarrow$
is a~unique correspondence from~$X$
to~$Y{\downarrow}$
satisfying the equality
$
\Psi\downwardarrow(x)=\Psi(x^{\scriptscriptstyle\wedge}){\downarrow}\quad (x\in X).
$
If
$\Phi:=f$
and~$\Psi:=g$
are functions then these relations take the
form
$
[\![f\upwardarrow(x^{\scriptscriptstyle\wedge})=f(x)]\!]={\mathbb 1}$
and $g\downwardarrow(x)=g(x^{\scriptscriptstyle\wedge})$ for all $x\in X$.

\subsection{2.7}
%Functional Representation of  Boolean Valued Universes}
Various function spaces reside in functional analysis,
and so the~problem is natural of replacing
an~abstract Boolean valued system by some function-space analog,
a~model whose elements are functions and in which the~basic logical operations
are calculated ``pointwise.''
An~example of such a~model is given by
the~class~${\mathbbm V}^Q$
of all functions defined on a~fixed nonempty
set~$Q$
and acting into~${\mathbbm V}$. The truth values on~${\mathbbm V}^Q$
are various subsets of~$Q$: The~truth value
$[\![\varphi(u_1,\dots,u_n)]\!]$
of
$\varphi(t_1,\dots,t_n)$
at functions
$u_1,\dots,u_n\in {\mathbbm V}^Q$
is calculated as follows:
$$
[\![\varphi(u_1,\dots,u_n)]\!]= \big\{q\in Q :\
 \varphi\big(u_1(q),\dots,u_n(q)\big)\big\}.
 $$

A.~G. Gutman and G.~A. Losenkov solved  the~above problem
by the concept
of continuous polyverse
which is a~continuous bundle of models of set theory.
It~is shown that the~class of continuous sections of a~continuous polyverse
is a~Boolean valued system satisfying all basic principles
of Boolean valued analysis and, conversely, each
Boolean valued algebraic system can be represented as the~class
of sections of a~suitable continuous polyverse.
More details are collected in~\cite[Chapter~6]{IBA}.

\ssection{3. Boolean Valued Algebraic Systems}

Every Boolean valued universe has  the collection of mathematical
objects in full supply:  available in plenty are all sets with
extra structure: groups, rings, algebras, normed spaces, etc.
Applying the descent functor to such {\it internal\/} algebraic
systems of a~Boolean valued model, we distinguish some bizarre
entities or recognize old acquaintances, which leads  to
revealing the new facts of their life and structure.

This technique of research, known as {\it direct Boolean valued
interpretation}, allows us to produce new theorems or, to be
more exact, to extend the semantical content of the available
theorems by means of  slavish translation.  The information we
so acquire might fail to be vital, valuable, or intriguing, in
which case the direct Boolean valued interpretation  turns out
into a leisurely game.

It thus stands to reason to raise  the following questions:
What structures significant for mathematical practice are
obtainable by the Boolean valued interpretation of the most
typical algebraic systems?  What transfer principles hold true
in this process?  Clearly, the answers should imply specific
objects whose particular  features enable us to deal with their
Boolean valued representation which, if understood duly, is
impossible to implement for arbitrary algebraic systems.

\subsection{3.1}~
An~{\it abstract Boolean set\/}
or~{\it set with ${\mathbbm B}$-structure\/}
is a~pair
$(X,d)$,
where
$X\in {\mathbbm V}$,
$X\ne \varnothing$,
and
$d$
is a~mapping from~$X\times X$
to~${\mathbbm B}$
such that
$d (x,y)={{\mathbb 0}}\leftrightarrow x=y$;
\
$d (x,y)=d (y,x)$;
\
$d (x,y)\le d (x,z)\vee d (z,y)$
all $x,y,z\in X$.

To obtain an easy~example of an abstract ${\mathbbm B}$-set, given
$\varnothing\ne X\subset {\mathbbm V}^{({\mathbbm B})}$
put
$$
d (x,y):=[\![x\ne y]\!]=\neg [\![x=y]\!]
$$ for $x,\ y\in X$.

Another easy example is a~nonempty $X$ with the
{\it discrete ${\mathbbm B}$-metric} $d$; i.e.,
$d (x,y)={\mathbb 1}$ if $x\ne y$ and $d (x,y)={\mathbb 0}$ if $x=y$.

\subsection{3.2}~Let
$(X,d)$
be some abstract ${\mathbbm B}$-set. There exist an~element
${\mathscr X}\in {\mathbbm V}^{({\mathbbm B})}$
and an~injection~$\iota:X\to X':={\mathscr X}{\downarrow}$
such that
$d (x,y)=[\![\iota x\ne \iota y]\!]$
for all $x, y\in X$
and every element~$x'\in X'$
admits the representation
$x'=\mix_{\xi\in \Xi}(b_{\xi}\iota x_{\xi})$,
where
$(x_{\xi})_{\xi\in \Xi}\subset X$
and
$(b_{\xi})_{\xi\in \Xi}$
is a~partition of unity in~${\mathbbm B}$.

We see that an abstract ${\mathbbm B}$-set
$X$  embeds in the~Boolean valued universe ${\mathbbm V}^{({\mathbbm B})}$ so
that the Boolean distance between the members of $X$ becomes the
Boolean truth value of the negation of their equality.  The
corresponding element ${\mathscr X}\in{\mathbbm V}^{({\mathbbm B})}$
is, by definition, the {\it Boolean valued representation\/} of~$X$.

If $X$ is a~discrete abstract ${\mathbbm B}$-set then
${\mathscr X}=X^{\scriptscriptstyle\wedge}$ and
$\iota x=x^{\scriptscriptstyle\wedge}$ for all $x\in X$.
If $X\subset {\mathbbm V}^{({\mathbbm B})}$ then $\iota{\uparrow}$
is an~injection from~$X{\uparrow}$ to~${\mathscr X}$
(inside~${\mathbbm V}^{({\mathbbm B})}$).

\subsection{3.3}~A~mapping~$f$
from a~${\mathbbm B}$-set
$(X,d)$
to a~${\mathbbm B}$-set
$(X',d')$
is said to be {\it contractive\/}
if
$d (x,y)\ge d'(f(x),f(y))$
for all
$x,y\in X$.

Let
$X$
and
$Y$
be some ${\mathbbm B}$-sets,
$\mathscr X$
and
$\mathscr Y$
be their Boolean-value representations, and
$\iota$
and
$\varkappa$
be the corresponding injections
$X\to \mathscr X{\downarrow}$
and
$Y\to \mathscr Y{\downarrow}$.
If
$f:X\to Y$
is a~contractive mapping then there is a~unique
element~$g\in{\mathbbm V}^{(\mathbbm B)}$
such that
$[\![g:\mathscr X\to \mathscr Y]\!]=\mathbb 1$
and
$f=\varkappa^{-1}\circ g{\downarrow}\circ\iota$.
We also accept the notations
$\mathscr X\!\!:=\mathscr F^\sim(X)\!\!:=X^\sim$
and
$g\!\!:=\mathscr F^\sim(f)\!\!:=f^\sim$.

\subsection{3.4}~The following are valid:

\subsubsec{(1)}~${\mathbbm V}^{({\mathbbm B})}\models f{(A)^\sim}={f^\sim}({A^\sim} )$
for $A\subset X$;

\subsubsec{(2)}~If $g:Y\to Z$ is a~contraction then
$g\circ f$ is a~contraction  and
${\mathbbm V}^{({\mathbbm B})}\models(g\circ f)^\sim ={g^\sim}\circ{f^\sim} $;

\subsubsec{(3)}~${\mathbbm V}^{({\mathbbm B})}\models$
``${f^\sim}$ is injective'' if and only if $f$ is
a~${\mathbbm B}$-isometry;

\subsubsec{(4)}~${\mathbbm V}^{({\mathbbm B})}\models$
``${f^\sim} $ is surjective'' if and only if
$\bigvee \{d(f(x),y):x\in X\}={\mathbbm 1}$
for every $y\in Y$.

\subsection{3.5}~In case a~${\mathbbm B}$-set $X$ has
some  a~priori structure we may try to furnish the Boolean
valued representation of $X$ with an analogous structure, so as
to apply the technique of ascending and descending to the study
of the original structure of~$X$.  Consequently, the above
questions may be treated as instances of the~unique problem of
searching a~well-qualified Boolean valued representation of
a~${\mathbbm B}$-set  with some additional structure.
We call these objects {\it algebraic ${\mathbbm B}$-systems\/}.

Recall that a~{\it signature\/}
 is a~3-tuple $\sigma\!\!:=(F,P,\mathfrak a)$, where
$F$ and $P$ are some (possibly, empty) sets and
$\mathfrak a$ is a~mapping from $F\cup P$ to~ $\omega$.
If the sets $F$ and $P$ are finite then $\sigma$
is a {\it finite signature}.
In applications we usually deal with algebraic systems
of finite signature.

An $n$-{\it ary operation\/}
and
an {\it $n$-ary predicate\/}
on a~$\mathbbm B$-set
$A$ are contractive mappings $f:A^n\to A$
and $p:A^n\to \mathbbm B$ respectively.
By definition,  $f$ and $p$
{\it are contractive mappings\/}
provided that
$$
\allowdisplaybreaks
\gathered
d(f(a_0,\dots,a_{n-1}),f(a'_0,\dots,a'_{n-1}))\le
\bigvee\limits_{k=0}^{n-1}\,d(a_k,a'_k),
\\
d_s\big(p(a_0,\dots,a_{n-1}),
p(a'_0,\dots,a'_{n-1})\big)
\le\bigvee\limits_{k=0}^{n-1}\,
d(a_k,a'_k)
\endgathered
$$
for all $a_0$, $a'_0,\dots,a_{n-1}$,
$a'_{n-1}\in A$, where $d$ is the  $\mathbbm B$-metric of  $A$, and $d_s$
is the {\it symmetric difference\/}
on~$\mathbbm B$; i.e.,
$d_s(b_1,b_2)\!\!:=b_1\triangle b_2 $
(cf. 1.1.4).

Clearly, the above definitions depend on
$\mathbbm B$ and it would be cleaner to speak of $\mathbbm B$-operations, $\mathbbm B$-predicates, etc.
We adhere to a simpler practice whenever it  entails no confusion.

\subsection{3.6}~An {\it algebraic $\mathbbm B$-system\/} $\mathfrak A$
of signature
$\sigma$ is a~pair
$(A,\nu)$,
where $A$ is a~nonemp\-ty
$\mathbbm B$-set,  the
{\it underlying set},
or {\it carrier}, or {\it universe\/}
of $\mathfrak A$,
and $\nu$ is a~mapping such that
(a)~$\dom(\nu)=F\cup P$;
(b) $\nu (f)$ is an $\mathfrak a(f)$-ary operation on ~$A$ for all
$f\in F$; and
(c)~ $\nu (p)$ is an $\mathfrak a(p)$-ary predicate on ~$A$
for every $p\in P$.

It is in common parlance to call $\nu$
the {\it interpretation\/}
of $\mathfrak A$, in which case the notation  $f^\nu$ and $p^\nu$ are common
substitutes for
$\nu (f)$ and $\nu (p)$.

The signature of an algebraic $\mathbbm B$-system $\mathfrak A\!\!:=(A,\nu)$
is often denoted by $\sigma(\mathfrak A)$;
while the carrier $A$ of $\mathfrak A$, by $|\mathfrak A|$.
Since  $A^0 =\{\varnothing\}$,  the nullary operations
and predicates on
$A$ are mappings from
$\{\varnothing\}$ to~$A$ and~$\mathbbm B$
respectively.
We agree to identify a mapping $g:\{\varnothing\}\to A\cup \mathbbm B$
with the element $g(\varnothing)$.
Each nullary operation on $A$ thus transforms into a unique
member of $A$. Analogously, the~set of all nullary predicates on
$A$ turns into the ~Boolean algebra ~$\mathbbm B$.
If $F\!\!:=\{f_1,\dots,f_n\}$ and
$P\!\!:=\{p_1,\dots,p_m\}$ then an  algebraic $\mathbbm B$-system
of signature $\sigma$ is often written down as
$(A,\nu (f_1),\dots,\nu (f_n)$,
$\nu (p_1),\dots,\nu (p_m))$
or even
$(A,f_1,\dots,f_n$, $p_1,\dots,p_m)$.
In this event, the expression $\sigma=(f_1,\dots,f_n$,
$p_1,\dots,p_m)$
is substituted for
$\sigma=(F,P,\mathfrak a)$.

\subsection{3.7}~We now address the $\mathbbm B$-valued interpretation
of a  first-order language. Consider an algebraic $\mathbbm B$-system
$\mathfrak A\!\!:=(A,\nu)$ of signature
$\sigma\!\!:=\sigma (\mathfrak A)\!\!:=(F,P,\mathfrak a)$.
Let $\varphi(x_0,\dots,x_{n-1})$ be a~formula of signature
$\sigma$ with
$n$ free variables. Assume given
$a_0,\dots,a_{n-1}\in A$.
We may readily define the truth value
$|\varphi |^{\mathfrak A}
(a_0,\dots,a_{n-1})\in \mathbbm B$
of a formula $\varphi$ in the system $\mathfrak A$ for the given values
$a_0,\dots,a_{n-1}$ of the variables $x_0,\dots,x_{n-1}$.
The definition proceeds as usual by induction on
the  complexity of~ $\varphi$:
Considering propositional connectives and quantifiers, we put
$$
\allowdisplaybreaks
\gathered
|\varphi\wedge\psi |^\mathfrak A\, (a_0,\dots,a_{n-1})\!:=
|\varphi |^\mathfrak A (a_0,\dots,a_{n-1})\wedge |\psi |^\mathfrak A
(a_0,\dots,a_{n-1});
\\
|\varphi\vee\psi |^\mathfrak A\,(a_0,\dots,a_{n-1})\!:=
|\varphi |^\mathfrak A (a_0,\dots,a_{n-1})\vee |\psi |^\mathfrak A
(a_0,\dots,a_{n-1});
\\
|\neg\varphi |^\mathfrak A\,(a_0,\dots,a_{n-1})\!:=|\varphi |^\mathfrak A (a_0,
\dots,a_{n-1})^* ;
\\
|(\forall\,x_0)\varphi |^\mathfrak A\, (a_1,\dots,a_{n-1})\!:=
\bigwedge\limits_{a_0\in A}|\varphi |^\mathfrak A (a_0,\dots,a_{n-1});
\\
|(\exists\,x_0)\varphi |^\mathfrak A\,(a_1,\dots,a_{n-1})\!:=
\bigvee\limits_{a_0\in A}|\varphi |^\mathfrak A (a_0,\dots,a_{n-1}).
\endgathered
$$
Now, the case of atomic formulas is in order.
Assume that  $p\in P$ symbolizes an $m$-ary predicate,
$q\in P$ is a~nullary predicate, and
$t_0,\dots,t_{m-1}$ be terms of signature
$\sigma$ assuming values $b_0,\dots,b_{m-1}$
at the given values $a_0,\dots,a_{n-1}$
of the variables $x_0,\dots,x_{n-1}$.
By definition, we let
$$
\gathered
|\varphi |^\mathfrak A (a_0,\dots,a_{n-1})\!:=\nu (q),
\text{~if~}\varphi=q^\nu;
\\
|\varphi |^\mathfrak A (a_0,\dots,a_{n-1})\!:=d(b_0,b_1)^*,
\,\kern5pt\text{~if~}\varphi=(t_0 =t_1);
\\
|\varphi |^\mathfrak A (a_0,\dots,a_{n-1})\!:=p^\nu
(b_0,\dots,b_{m-1}),
\,\kern5pt\text{~if~}\varphi=p^\nu (t_0,\dots,t_{m-1}),
\endgathered
$$
where $d$ is a~ $\mathbbm B$-metric on  $A$.

Say that
$\varphi(x_0,\dots,x_{n-1})$
is {\it valid\/}
in
$\mathfrak A$ at the given  values
$a_0,\dots,a_{n-1}\in A$
of  $x_0,\dots,x_{n-1}$
and write $\mathfrak A\models\varphi(a_0,\dots,a_{n-1})$
provided that
$|\varphi |^\mathfrak A (a_0,\dots,a_{n-1})=\mathbb 1_\mathbbm B$.
Alternative expressions
are as follows:
$a_0,\dots,a_{n-1}\in A$
{\it satisfies\/}
$\varphi(x_0,\dots,x_{n-1})$;
or $\varphi(a_0,\dots,a_{n-1})$
holds true in
$\mathfrak A$.
In case $\mathbbm B\!:=\{\mathbb 0,\mathbb 1\}$, we arrive at the conventional
definition of
the validity of a~formula in an  algebraic system.
%(cf. [\cite{33}, \cite{106}]).

Recall that a~closed formula
$\varphi $
of signature $\sigma$ is
{\it tautology\/}
if
$\varphi$ is valid on every algebraic
$\mathbb 2$-system of signature $\sigma$.

\subsection{3.8}~Before giving a~general definition of the descent of an
algebraic system, consider the descent
of a~very simple but important algebraic system,
the two-element Boolean algebra.
Choose two arbitrary elements,
$0$, $1\in\mathbbm V ^{(\mathbbm B)}$,
satisfying $[\![0\ne 1]\!]=\mathbb 1_{\mathbbm B}$.
We may for instance assume that
$0\!:=\mathbb 0_{\mathbbm B}^{\scriptscriptstyle\wedge}$
and
$1\!:=\mathbb 1^{\scriptscriptstyle\wedge} _\mathbbm B$.

\Proclaim{}
The descent $C$ of the two-element Boolean algebra
$\{0,1\}^{\mathbbm B}\in\mathbbm V ^{(\mathbbm B)}$
is a~complete Boolean algebra isomorphic to  ~$\mathbbm B$.
The formulas
$$
[\![\chi(b)=1]\!]=b,\quad
[\![\chi (b)=0]\!]=b^*\quad
(b\in \mathbbm B)
$$
defines an isomorphism $\chi:\mathbbm B\to C$.
\endproc

\subsection{3.9}~Consider now an algebraic system $\mathfrak A$ of signature
$\sigma^{\scriptscriptstyle\wedge}$ inside $\mathbbm V ^{(\mathbbm B)}$, and let
$[\![\mathfrak A =(A,\nu)^\mathbbm B]\!]=\mathbb 1$ for some
$A$, $\nu\in\mathbbm V ^{(\mathbbm B)}$.
The {\it descent\/} of
$\mathfrak A$ is the pair
$\mathfrak A\!\!\downarrow\!:=(A\!\!\downarrow,\mu)$, where $\mu $
is the~function  determined from the formulas:
$$
\gathered
\mu:f\mapsto(\nu\downwardarrow (f))\!\!\downarrow\quad
(f\in F),
\\
\mu:p\mapsto\chi ^{-1}\circ (\nu\downwardarrow (p))\!\!\downarrow\quad
(p\in P).
\endgathered
$$
Here $\chi$ is the above-defined isomorphism of~$\mathbbm B$.

In more detail, the modified descent
$\nu\downwardarrow$ is the~mapping with domain
$\dom(\nu\downwardarrow)=F\cup P$.
Given $p\in P$, observe $[\![\mathfrak a (p)^{\scriptscriptstyle\wedge}
=\mathfrak a^{\scriptscriptstyle\wedge} (p^{\scriptscriptstyle\wedge})]\!]=\mathbb 1$,
$[\![\nu\downwardarrow (p)=\nu (p^{\scriptscriptstyle\wedge})]\!]=\mathbb 1$
and so
$$
\mathbbm V ^{(\mathbbm B)}\models\nu\downwardarrow (p):A^{\mathfrak a
(f)^{\scriptscriptstyle\wedge}}
\to\{0,1\}^\mathbbm B.
$$
It is now obvious that $(\nu\downwardarrow (p))\!\!\downarrow:
(A\!\!\downarrow)^{\mathfrak a (f)}
\to C\!:=\{0,1\}^\mathbbm B\!\!\downarrow $
and we may put
$\mu (p)\!:=\chi ^{-1}\circ (\nu\downwardarrow (p))\!\!\downarrow$.

\subsection{3.10}~Let $\varphi(x_0,\dots,x_{n-1})$ be a~fixed formula of signature
$\sigma$
in $n$ free variables.
Write down the~formula $\Phi (x_0,\dots,x_{n-1},\!\mathfrak A)$
in the language of set theory which formalizes the proposition
$\mathfrak A\models\varphi(x_0$,
$\dots,x_{n-1})$.
Recall that the  formula
$\mathfrak A\models\varphi(x_0,\dots,x_{n-1})$
determines an $n$-ary predicate on~$A$ or, which is the same, a~mapping
from $A^n$ to $\{0,1\}$.
By the  maximum and transfer principles, there is a~unique element
$|\varphi |^\mathfrak A\in\mathbbm V ^{(\mathbbm B)}$
such that
$$
%\gathered
[\![|\varphi |^\mathfrak A:A^{n^{\scriptscriptstyle\wedge}}\to\{0,1\}^{\mathbbm B}]\!]=\mathbb 1,
$$
$$
%\\
[\![|\varphi |^\mathfrak A (a\!\!\uparrow)=1]\!]=
[\![\Phi (a(0),\dots,a(n-1),\mathfrak A)]\!]=\mathbb 1
%\endgathered
$$
for every $a:n\to A\!\!\downarrow$.
Henceforth instead of $|\varphi |^\mathfrak A$ $(a\!\!\uparrow)$
we will write
$|\varphi |^\mathfrak A$ $(a_0,\dots,a_{n-1})$, where
$a_l\!:=a(l)$.
Therefore, the formula
$$
\mathbbm V ^{(\mathbbm B)}\models\text{``}\varphi(a_0,\dots,a_{n-1})
\text{~is valid in~}\mathfrak A\text{''}
$$
holds true if and only if
$[\![\Phi (a_0,\dots,a_{n-1},\mathfrak A)]\!]=\mathbb 1$.

\Proclaim{}
Let $\mathfrak A$ be an algebraic system of signature
$\sigma^{\scriptscriptstyle\wedge} $
inside $\mathbbm V ^{(\mathbbm B)}$.
Then $\mathfrak A\!\!\downarrow$ is a ~universally complete algebraic
$\mathbbm B$-system of signature ~$\sigma$.
In this event,
$$
\chi\circ |\varphi |^{\mathfrak A\downarrow}=|\varphi |^\mathfrak A\!\!\downarrow.
$$
foe each formula $\varphi$ of signature
$\sigma$.
\endproc

\subsection{3.11}
{\sl
Let  $\mathfrak A\!:=(A,\nu)$ be an algebraic $\mathbbm B$-system of signature $\sigma$.
Then there are  $\mathscr A$ and $\mu\in\mathbbm V ^{(\mathbbm B)}$ such that
the following  are fulfilled:

\subsubsec{(1)}
$\mathbbm V ^{(\mathbbm B)}\models $
``$(\mathscr A,\mu)$ is an algebraic system of signature
$\sigma^{\scriptscriptstyle\wedge} $'';

\subsubsec{(2)}
If $\mathfrak A '\!:=(A',\nu ')$ is the descent of
$(\mathscr A,\mu)$ then $\mathfrak A '$ is a ~universally complete algebraic
$\mathbbm B$-system of signature $\sigma$;

\subsubsec{(3)}
There is an isomorphism $\imath$ from $\mathfrak A$ to $\mathfrak A '$ such that
$A'=\mix(\imath(A))$;

\subsubsec{(4)}
For every formula $\varphi$ of signature $\sigma$ in  $n$ free variables,
the equalities hold
$$
\gathered
|\varphi |^\mathfrak A (a_0,\dots,a_{n-1})=|\varphi |^{\mathfrak A '}
(\imath(a_0),\dots,\imath(a_{n-1}))
\\
=\chi ^{-1}\circ
(|\varphi |^{\mathfrak A^\sim})\!\!\downarrow\!\!(\imath(a_0),\dots,
\imath(a_{n-1}))
\endgathered
$$
for all $a_0,\dots,a_{n-1}\in A$ and $\chi$ the same as in~3.8.
}

\subsection{3.12}
The Boolean valued representation of an algebraic ${\mathbbm B}$-system
appears to be a~conventional two-valued algebraic system of the
same type.  This means that an appropriate completion of each
algebraic ${\mathbbm B}$-system coincides with the descent of
some two-valued algebraic system inside ${\mathbbm V} ^{({\mathbbm B})}$.

On the other hand, each two-valued algebraic system  may be
transformed into an algebraic ${\mathbbm B}$-system on distinguishing
a~complete Boolean algebra of congruences of the original
system.  In this event, the task is in order of finding the
formulas  holding true in direct or reverse transition from
a~${\mathbbm B}$-system to a~two-valued system.  In other words, we  have
to seek here for some versions of the transfer  or
identity preservation principle of long standing in some
branches of mathematics.

%\newpage
\ssection{4. Boolean Valued Numbers}

Boolean valued analysis stems from the fact that each internal
field of reals of a~Boolean valued model descends into
a~universally complete Kantorovich space.  Thus, a~remarkable
opportunity opens up to expand and enrich the treasure-trove of
mathematical knowledge by translating information about the
reals to the language of other noble families of functional
analysis. We will elaborate  upon the matter in  this section.

\subsection{4.1}~Recall a~few definitions. A {\it vector
lattice\/} is an~ordered vector space whose order makes it a~lattice.
In other words,  the {\it join\/} $\sup \{x_1, \dots, x_n\}\!:= x_1 \vee \dots \vee x_n$ and
{\it meet\/} $\inf \{x_1, \dots, x_n\}\!:= x_1 \wedge
\dots \wedge x_n$ correspond to each finite subset
$\{x_1, \dots, x_n\}$ of~a vector lattice.  In particular, each
element~$x$ has the {\it positive part}~ $x^+ \!:= x\vee 0$,
{\it negative part\/} $x^- \!:=(-x)^+ \!:= -x \wedge 0$, and
{\it modulus}~$|x| \!:= x\vee (-x)$.

A vector lattice $E$ is called {\it Archimedean\/} if for every
pair of elements $x,y\in E$ from $(\forall n\in \mathbb N)\
nx\le y$ it follows that $x\le 0$.
We  assume all vector lattices  Archimedean in what follows.

Two elements $x$ and $y$ of a~vector lattice $E$ are  {\it
disjoint\/} (in symbols $x\perp y$) if $|x| \wedge  |y| =0$.
A~{\it band\/} of $E$ is defined as the {\it disjoint
complement\/} $M^\perp\!:= \{x\in E :\, (\forall y\in M)\,
x\perp y\}$ of a~nonempty set $M\subset E$.

The inclusion-ordered set ${\mathfrak B}(E)$ of all bands in $E$
is a~complete Boolean algebra with the Boolean operations:
$$
L\wedge K=L\cap K,\quad
L\vee K=(L\cup K)^{\perp\perp},\quad
L^* =L^\perp\quad
(L,K\in{\mathfrak B}(E)).
$$
The Boolean algebra ${\mathfrak B}(E)$ is often referred as to
the {\it base\/} of~$E$.

A~{\it band projection\/} in $E$ is a~linear idempotent operator
in $\pi:E\to E$ satisfying the inequalities $0\leq\pi x\leq x$
for all $0\leq x\in E$. The set ${\mathfrak P}(E)$ of all band
projections ordered by $\pi\le\rho\Longleftrightarrow\pi\circ\rho=\pi$
is a~Boolean algebra with the Boolean operations:
$$
\pi\wedge \rho =\pi\circ\rho,\quad
\pi\vee\rho =\pi +\rho -\pi\circ\rho,\quad
\pi^* =I_E-\pi\quad (\pi,\rho\in\mathfrak (E)).
$$

Let $u\in E_+$ and $e\wedge (u-e)= 0$ for some $0\leq e\in E$.
Then $e$ is  a~{\it fragment\/} or {\it component\/} of~$u$.
The set ${\mathfrak E}(u)$ of all fragments of $u$ with the order
induced by~$E$ is a~Boolean algebra where the lattice operations
are taken from~$E$ and the Boolean complement has the form
$e^* \!:= u - e$.

\subsection{4.2}~A Dedekind complete vector lattice is also
called a~{\it Kantorovich space\/} or $K$-{\it space\/}, for
short.  A $K$-space $E$ is {\it universally complete\/} if
every family of pairwise disjoint elements of~$E$ is order
bounded.

\subsubsec{(1)}
\proclaim{Theorem.} Let $E$ be an arbitrary \hbox{$K$-}space. Then the correspondence
$\pi\mapsto\pi(E)$ determines an~isomorphism of the Boolean
algebras~${\mathfrak P}(E)$ and~${\mathfrak B}(E)$. If there is an~order unity
${\mathbb 1}$ in~$E$ then the mappings $\pi\mapsto\pi{\mathbb 1}$
from~${\mathfrak P}(E)$ into~${\mathfrak E}(E)$ and
$e\mapsto \{e\}^{\perp\perp}$ from~${\mathfrak E}(E)$ into~${\mathfrak B}(E)$
are isomorphisms of Boolean algebras too.
\endproc

\subsubsec{(2)}
\proclaim{Theorem.}
Each universally complete $K$-space $E$ with order unity~${\mathbb 1}$
can be uniquely endowed by multiplication so as to make~$E$
into a~faithful $f$-algebra and ${\mathbb 1}$ into a~ring unity.
In this $f$-algebra each band projection $\pi\in{\mathfrak P}(E)$
is the operator of multiplication by~$\pi({\mathbb 1})$.
\endproc

\subsection{4.3}
By a~{\it field of reals\/} we mean every algebraic
system that satisfies the axioms of an~Archimedean  ordered
field (with distinct zero and unity) and enjoys the axiom of
completeness. The same object can be defined as a~one-dimensional
$K$-space.

Recall the well-known assertion of $\ZFC$:
{\sl There exists a~field of reals ${\mathbbm R}$ that is unique
up to isomorphism.}

Successively applying the transfer and maximum principles,
we find an~element
${\mathscr R}\in{\mathbbm V}^{({\mathbbm B})}$ for which $[\![\,{\mathscr R}$ is a~field of
reals$\,]\!]={\mathbb 1}$.  Moreover, if an~arbitrary
${\mathscr R}\,'\in {\mathbbm V}^{({\mathbbm B})}$ satisfies the condition
$[\![{\mathscr R}\,'$ is a~field of  reals$\,]\!]={\mathbb 1}$ then
$[\![\,$the ordered fields $\mathscr R$ and ${\mathscr R}\,'$
are isomorphic\,$]\!]=\nobreak{\mathbb 1}$.  In other words, there exists an
internal field~of  reals ${\mathscr R}\in{\mathbbm V}^{({\mathbbm B})}$ which
is unique up to isomorphism.

By the same reasons there exists an internal field~of complex
numbers ${\mathscr C}\in{\mathbbm V}^{({\mathbbm B})}$ which is
unique up to isomorphism. Moreover,
${\mathbbm V}^{({\mathbbm B})}\models{\mathscr C}={\mathscr R}\oplus i{\mathscr R}$.
We call ${\mathscr R}$ and ${\mathscr C}$ the {\it internal  reals} and {\it internal complexes\/}
in ${\mathbbm V}^{({\mathbbm B})}$.

\subsection{4.4}~Consider another well-known assertion of
$\ZFC$:  {\sl If ${\mathbbm P}$ is an Archimedean ordered field
then there is an~isomorphic embedding~$h$ of the
field~${\mathbbm P}$ into~${\mathbbm R}$ such that the
image~$h({\mathbbm P})$ is a~subfield of~${\mathbbm R}$
containing the subfield of rational numbers.  In particular,
$h({\mathbbm P})$ is dense in~${\mathbbm R}$.}

Note also that $\varphi (x)$, presenting the conjunction
of the axioms of an~Archimedean ordered field $x$, is bounded;
therefore, $[\![\,\varphi ({\mathbbm R}^{\scriptscriptstyle\wedge})\,]\!]
={\mathbb 1}$,
i.e., $[\![\,{\mathbbm R}^{\scriptscriptstyle\wedge}$ is an~Archimedean ordered
field$\,]\!]={\mathbb 1}$. ``Pulling'' 4.2\,(2) through
the transfer principle, we conclude that
$[\![\,{\mathbbm R}^{\scriptscriptstyle\wedge}$ is isomorphic to
a~dense subfield of~$\mathscr R\,]\!]={\mathbb 1}$. We further
assume that ${\mathbbm R}^{\scriptscriptstyle\wedge}$ is a~dense
subfield of ${\mathscr R}$ and ${\mathbbm C}^{\scriptscriptstyle\wedge}$
is a~dense subfield of ${\mathscr C}$. It is easy to note that the
elements~$0^{\scriptscriptstyle\wedge}$
and~$1^{\scriptscriptstyle\wedge}$ are the zero and
unity of~${\mathscr R}$.

Observe that the equalities
${\mathscr R}={\mathbbm R}^{\scriptscriptstyle\wedge}$ and
${\mathscr C}={\mathbbm C}^{\scriptscriptstyle\wedge}$ are not valid in
general. Indeed, the axiom of completeness for~${\mathbbm R}$ is
not a~bounded formula and so it may thus fail for
${\mathbbm R}^{\scriptscriptstyle\wedge}$
inside~${\mathbbm V}^{({\mathbbm B})}$.

\subsection{4.5}~Look now at the descent
${\mathscr R}{\downarrow}$
of the algebraic system~${\mathscr R}$.
In other words,  consider the descent of the underlying set of the
system~${\mathscr R}$ together with descended operations and order. For
simplicity, we denote the operations and order in~${\mathscr R}$
and~${\mathscr R}{\downarrow}$ by the same symbols $+$, $\cdot\,$, and $\le $.
In more detail, we introduce addition, multiplication, and order
in~${\mathscr R}{\downarrow}$
by the formulas
$$
\gathered
z=x+y\leftrightarrow [\![\,z=x+y\,]\!]={\mathbb 1},
\\
z=x\cdot y\leftrightarrow [\![\,z=x\cdot y\,]\!]={\mathbb 1},
\\
x\le y\leftrightarrow [\![\,x\le y\,]\!]={\mathbb 1}
\quad
(x, y, z\in {\mathscr R}{\downarrow}).
\endgathered
$$
Also, we may introduce multiplication by the usual reals
in~${\mathscr R}{\downarrow}$ by the rule
$$
y=\lambda x\leftrightarrow
[\![\,\lambda^{\scriptscriptstyle\wedge}x=y\,]\!]={\mathbb 1}\quad
(\lambda\in{\mathbbm R},\ x, y\in {\mathscr R}{\downarrow}).
$$

One on the most fundamental results of Boolean valued analysis
reads: {\sl Each universally complete
Kantorovich space is an~interpretation of the reals in
an~appropriate Boolean valued model}. In other words, we have the following

\subsection{4.6}
\proclaim{Gordon Theorem.}
Let
${\mathscr R}$
be the reals inside~${\mathbbm V}^{({\mathbbm B})}$.
Then
${\mathscr R}{\downarrow}$,
(with the descended operations and order, is a~universally complete
$K$-space with order unity~$1$.
Moreover, there exists an~isomorphism~$\chi$
of~${\mathbbm B}$
onto~${\mathfrak P} ({\mathscr R}{\downarrow})$
such that
$$
\chi (b) x=\chi (b) y\leftrightarrow b\le [\![\,x=y\,]\!],
\quad
\chi (b) x\le \chi (b) y\leftrightarrow b\le [\![\,x\le y\,]\!]
$$
for all
$x, y\in {\mathscr R}{\downarrow}$
and
$b\in {\mathbbm B}$.
\endproc

The converse  is also true: {\sl Each Archimedean vector lattice
embeds in a~Boolean valued model, becoming a~vector sublattice
of the reals} (viewed as such over some dense subfield of the reals).

\subsection{4.7}
\proclaim{Theorem.}
Let
$E$
be an~Archimedean vector lattice, let
${\mathscr R}$
be the reals inside~${\mathbbm V}^{({\mathbbm B})}$,
and let
$\jmath$
be an isomorphism of~${\mathbbm B}$
onto~${\mathfrak B}(E)$.
Then there is ${\mathscr E}\in {\mathbbm V}^{({\mathbbm B})}$
such that

\subsubsec{(1)}
${\mathscr E}$
is a~vector sublattice of~${\mathscr R}$
over~${\mathbbm R}^{\scriptscriptstyle\wedge}$ inside ${\mathbbm V}^{({\mathbbm B})}$;

\subsubsec{(2)}
$E'\!:={\mathscr E}{\downarrow}$
is a~vector sublattice of~${\mathscr R}{\downarrow}$
invariant under every band projection
$\chi (b)$
$(b\in {\mathbbm B})$
and such that each set of positive pairwise disjoint
elements in it has a~supremum;

\subsubsec{(3)}
there is an~\hbox{$o$-}continuous
lattice isomorphism~$\iota:E\to E'$
such that $\iota (E)$
is a~coinitial sublattice of~${\mathscr R}{\downarrow}$;

\subsubsec{(4)}
for every~$b\in {\mathbbm B}$ the band projection
in~${\mathscr R}{\downarrow}$
onto
$\{\iota (\jmath(b))\}^{\perp\perp}$
coincides with~$\chi (b)$.
\endproc

Note also that ${\mathscr E}$ and ${\mathscr R}$ coincide if
and only if $E$ is Dedekind complete.  Thus,
each theorem about the reals within Zermelo--Fraenkel set theory
has an analog in an arbitrary Kantorovich space. Translation of
theorems is carried out by appropriate general functors of
Boolean valued analysis. In particular, the most important structural
properties of vector lattices such as the functional representation,
spectral theorem,  etc. are the ghosts of some properties of the reals in
an~appropriate Boolean valued model.
More details and references are collected in~\cite{IBA}.

\subsection{4.8}~The theory of vector lattices with a~vast
field of applications is thoroughly  covered in many monographs
(see~\cite{AK, AB,  KA, KVP, LZ1, Sch, Schw, Vul, Z}).
The credit for finding the most important instance among ordered vector
spaces, an~order complete vector lattice or~$K$-space, is due to
L.~V.~Kantorovich. This notion appeared in Kantorovich's first
article on this topic~\cite{Ka1} where he wrote: %\?
``In this note, I define a~new type of space that I call
a~semiordered linear space. The introduction of such a~space allows us to
study linear operations of one abstract class
(those with values in such a~space) as linear functionals.''

Thus the {\it heuristic transfer principle\/} was stated
for~$K$-spaces which becomes the Ariadna thread of many subsequent studies.
The depth and universality of Kantorovich's principle are explained within
Boolean valued analysis.

\subsection{4.9}~Applications of Boolean valued models to
functional analysis stem from the works by E.~I.~Gordon
\cite{Gor1, Gor2} and G.~Takeuti \cite{Tak}. If
${\mathbbm B}$ in~4.6 is the algebra of $\mu$-measurable sets modulo
$\mu$-negligible sets then ${\mathscr R}{\downarrow}$ is isomorphic
to the universally complete $K$-space~$L^0(\mu)$ of measurable
functions. This fact (for the Lebesgue measure on an~interval)
was already known to D.~Scott and R.~Solovay (see~\cite{IBA}).
If ${\mathbbm B}$ is a~complete Boolean algebra of projections in a~Hilbert
space then ${\mathscr R}{\downarrow}$ is isomorphic to the space of
selfadjoint operators~${\mathfrak A}({\mathbbm B})$.  These two
particular cases of Gordon's Theorem were intensively and
fruitfully exploited by G.~Takeuti (see~\cite{Tak} and the
bibliography in~\cite{IBA}). The object ${\mathscr R}{\downarrow}$
for general Boolean algebras was also studied by
T.~Jech~\cite{Jech1}--\cite{Jech3} who in fact rediscovered
Gordon's Theorem. The difference is that in~\cite{Jech}
a~(complex) universally complete $K$-space with unity is defined
by another system of axioms and is referred to as a~{\it complete
Stone algebra}. Theorem~4.7 was obtained
by A.~G.~Kusraev \cite{K9}. A~close result (in other terms)
is presented in T.~Jech's article~\cite{Jech3} where some
Boolean valued interpretation is revealed of the theory of linearly ordered
sets.  More details can be found in \cite{IBA}.

\ssection{5. Band Preserving Operators}

This section deals with the class of band preserving operators.
Simplicity of these operators notwithstanding, the~question
about their order boundedness is far from trivial.

\subsection{5.1}
Recall that a~complex $K$-space is the complexification
$G_{{\mathbbm C}}\!:=G\oplus iG$ of a~real $K$-space $G$ (see \cite{Sch}).
A~linear operator $T:G_{{\mathbbm C}}\to G_{{\mathbbm C}}$ is
{\it band preserving}, or {\it contractive}, or a~{\it stabilizer\/} if, for all $f,g\in G_{{\mathbbm C}}$, from
$f\perp g$ it follows that $Tf\perp g$. Disjointness in
$G_{{\mathbbm C}}$ is defined just as in $G$ (see 4.1), whereas
$|z|\!:=\sup\{\Re (e^{i\theta } z):\, 0\leq\theta\leq\pi\}$
for $z\in G_{{\mathbbm C}}$. Thus, a linear operator is band
preserving if every band is its invariant subspace.

\subsubsec{(1)}~Let $\End_N(G_{{\mathbbm C}})$ stand for the set
of all band preserving linear operators in $G_{{\mathbbm C}}$, with
$G\!:={\mathscr R}{\downarrow}$. Clearly, $\End_N(G_{{\mathbbm C}})$
is a~complex vector space. Moreover, $\End_N(G_{{\mathbbm C}})$
becomes a~faithful unitary module over the ring $G_{{\mathbbm C}}$
if we define $gT$ as $gT:x\mapsto g\cdot Tx$ for all $x\in G$.
This follows from the fact that multiplication by a~member of
$G_{{\mathbbm C}}$ is a~band preserving operator and the composite
of band preserving operators is band preserving too.

\subsubsec{(2)}~Denote by $\End_{{\mathbbm C}^{\scriptscriptstyle\wedge}}({\mathscr C})$
the element of~${\mathbbm V}^{({\mathbbm B})}$ representing the space of all
${\mathbbm C}^{\scriptscriptstyle\wedge}$-linear mappings from
${\mathscr C}$ to~${\mathscr C}$. Then
$\End_{{\mathbbm C}^{\scriptscriptstyle\wedge}}({\mathscr C})$
is a~vector space over ${\mathbbm C}^{\scriptscriptstyle\wedge}$
inside ${\mathbbm V}^{({\mathbbm B})}$, and
$\End_{{\mathbbm C}^{\scriptscriptstyle\wedge}}({\mathscr C}){\downarrow}$
is a~faithful unitary module over~$G_{{\mathbbm C}}$.

\subsection{5.2}~Following \cite{Kus1} it is easy to prove that
a linear operator $T$ in the $K$-space $G_{{\mathbbm C}}$
is band preserving if and only if $T$ is extensional.
Since  each extensional mapping has an ascent,
$T\in\End_N(G_{{\mathbbm C}})$ has the ascent $\tau\!:=T{\uparrow}$
which is a~unique internal functional from ${\mathscr C}$ to ${\mathscr C}$
such that
$[\![\tau (x)=Tx]\!]={\mathbb 1}$ $(x\in G_{{\mathbbm C}})$.
We thus arrive at the following assertion:

\Proclaim{}
The modules $\End_N(G_{{\mathbbm C}})$ of all linear
band preserving operators in the complex $K$-space $G_{\mathbbm C}$
and the descent of the internal space
$\End_{{\mathbbm C}^{\scriptscriptstyle\wedge}}({\mathscr C}){\downarrow}$
of ${\mathbbm C}^{\scriptscriptstyle\wedge}$-linear functions in the
internal complexes ${\mathscr C}$ (considered as a vector space over
${\mathbbm C}^{\scriptscriptstyle\wedge}$)
are isomorphic by sending each band preserving operator to its ascent.
\endproc

By Gordon's Theorem this assertion means that the problem of finding
a band preserving operator in $G_{{\mathbbm C}}$ amounts to solving
(for $\tau:{\mathscr C}\to{\mathscr C}$)
inside ${\mathbbm V}^{({\mathbbm B})}$ the {\it  Cauchy functional equation}:
$\tau(x+y)=\tau(x)+\tau (y)\quad(x,y\in{\mathscr C})$
under the subsidiary condition $\tau(\lambda x)=\lambda\tau (x)
\quad(x\in{\mathscr C},\lambda\in{\mathbbm C}^{\scriptscriptstyle\wedge})$.

As another subsidiary condition we may consider the {\it Leibniz rule}
$\tau(xy)=\tau(x)y+x\tau(y)$(in which case $\tau$ is called a~${\mathbbm
C}^{\scriptscriptstyle\wedge}$-{\it derivation}) or multiplicativity
$\tau(xy)=\tau(x)\tau(y)$. These situations are addressed in~5.5.

\subsection{5.3}~An~element $g\in G^{+}$ is  {\it locally constant\/}
with respect to~$f\in G^{+}$ if
$g=\bigvee _{\xi \in \Xi } \lambda _{\xi }\pi _{\xi }f$ for some
numeric family $(\lambda _{\xi })_{\xi \in \Xi }$ and a~family
$(\pi _{\xi })_{\xi \in \Xi }$ of pairwise disjoint band
projections. A~universally complete $K$-space~$G_{{\mathbbm C}}$
is called {\it locally one-dimensional\/} if all
elements of~$G^{+}$ are locally constant with respect to some
order unity of~$G$ (and hence each of them). Clearly,
a~$K$-space $G_{{\mathbbm C}}$ is locally one-dimensional
if each  $g\in G_{{\mathbbm C}}$ may be presented as
$g=\osum_{\xi \in \Xi}\lambda_{\xi }\pi_{\xi}{\mathbb 1}f$
with some family $(\lambda_{\xi })_{\xi \in \Xi}\subset{\mathbbm C}$
and partition of unity $(\pi _{\xi })_{\xi \in \Xi}\subset{\mathfrak P}(G)$.

\proclaim{Theorem}
Let $G_{{\mathbbm C}}$ be a~universally complete $K$-space. Every band
preserving linear operator in $G_{{\mathbbm C}}$ is order bounded
if and only if $G_{{\mathbbm C}}$ is locally one-dimensional.
\endproc

\subsection{5.4}~A~$\sigma$-complete Boolean
algebra~${\mathbbm B}$ is called {\it $\sigma$-distributive\/} if
$$
\bigvee_{n\in {{\mathbbm N}}}\bigwedge_{m\in {{\mathbbm N}}} b_{n,m}=
\bigwedge_{\varphi\in {{\mathbbm N}}^{{{\mathbbm N}}}}
\bigvee  _{n\in {{\mathbbm N}}} b_{n,\varphi(n)}.
$$
for every double sequence  $(b_{n,m})_{n,m\in{\mathbbm N}}$ in~${\mathbbm B}$.

An equivalent definition can be given in terms of partitions of unity.  From
any two partitions of unity in an arbitrary Boolean algebra one can refine a
partition of unity by taking infimum of any pair of members of the
partitions. The same is true for a finite set of partitions of unity. A
$\sigma $-complete Boolean algebra ${\mathbbm B}$ is called
$\sigma$-distributive if from every sequence of countable partitions of unity
in ${\mathbbm B}$, it is possible to refine a (possibly, uncountable)
partition of unity.

Other equivalent definitions are collected in~\cite{Sik}.
As an example of a~$\sigma$-distributi\-ve Boolean algebra
we may take a~complete atomic Boolean algebra, i.e., the boolean of
a nonempty set. It is worth observing that there are
nonatomic $\sigma$-distributive complete Boolean algebras
(see \cite[5.1.8]{DOP}).

\subsection{5.5}~We now address the problem which is often referred
to in the literature as {\it Wickstead's problem}: Characterize the
universally complete vector lattices in which every band preserving
linear operator is order bounded. We restrict exposition to the case
of complex vector lattices.

According to 5.2,  Boolean valued analysis reduces Wickstead's
problem to that of order boundedness of the endomorphisms
of the field ${\mathscr C}$ viewed as a~vector space and algebra
over~${\mathbbm C}^{\scriptscriptstyle\wedge}$. It is
important that the standard name of the external complexes is
an~algebraically closed field inside ${\mathbbm V}^{(\mathbbm B)}$:

\Proclaim{}
The field ${\mathbbm C}^{\scriptscriptstyle\wedge}$ is
algebraically closed in $\mathscr C$ inside ${\mathbbm V}^{(\mathbbm B)}$. In particular, if
${\mathbbm C}^{\scriptscriptstyle\wedge}\ne\mathscr C$ then
$$
{\mathbbm V}^{(\mathbbm B)}\models\text{``$\mathscr C$~---
transcendental extension of the field''~}{\mathbbm C}^{\scriptscriptstyle\wedge}.
$$
\endproc

We so arrived at an internal Cauchy type functional equation: Find an
additive function in the internal complexes that is ${\mathbbm
P}$-homogeneous for some algebraically closed dense subfield ${\mathbbm P}$.
(If ${\mathbbm P}$ coincides with the field of rationals then we obtain
exactly the Cauchy functional equation inside the Boolean valued universe.) The
corresponding scalar result reads as follows.

\subsection{5.6}
\proclaim{Theorem.}
Let ${\mathbbm P}$ be an~algebraically closed and (topologically) dense
subfield of the field of complexes ${\mathbbm C}$. The following
are equivalent:

\subsubsec{(1)}~${\mathbbm P}={\mathbbm C}$;

\subsubsec{(2)}~every ${\mathbbm P}$-linear function on ${\mathbbm C}$
is order bounded;

\subsubsec{(3)}~there are no nontrivial
${\mathbbm P}$-derivations on ${\mathbbm C}$;

\subsubsec{(4)}~each ${\mathbbm P}$-linear
endomorphism on ${\mathbbm C}$ is the zero or identity function;

\subsubsec{(5)}~there is no ${\mathbbm P}$-linear
automorphism on ${\mathbbm C}$ other than the identity.
\endproc

The equivalence (1)~$\leftrightarrow$~(2) is checked by
using a~Hamel basis of the vector space ${\mathbbm C}$ over~${\mathbbm P}$.
The remaining equivalences  rest on replacing a~Hamel basis with a
transcendence basis (for details see \cite{Kus2}).

Recall that a~linear operator $D:G_{{\mathbbm C}}\to G_{{\mathbbm C}}$
is a~${\mathbbm C}$-{\it derivation\/} if it obeys the Leibnitz rule
$D(fg)=D(f)g+fD(g)$ for all $f,g\in G_{{\mathbbm C}}$.
It can be easily checked that every ${\mathbbm C}$-derivation
is band preserving.

Interpreting  Theorem~5.5 in~${{\mathbbm V}}^{({\mathbbm B})}$,
we arrive at following two results.

\subsection{5.7}
\proclaim{Theorem.}
If ${\mathbbm B}$ is a~complete Boolean algebra then the
following are equivalent:

\subsubsec{(1)}~${\mathscr C}=
{{\mathbbm C}}^{\scriptscriptstyle\wedge}$ inside~${{\mathbbm V}}^{({\mathbbm B})}$;

\subsubsec{(2)}~every band preserving linear operator is order bounded in the complex
vector lattice~${\mathscr C}{\downarrow}$;

\subsubsec{(3)}~${\mathbbm B}$ is $\sigma$-distributive.
\endproc

\subsection{5.8}
\proclaim{Theorem.}
If ${\mathbbm B}$ is a~complete Boolean algebra then the
following are equivalent:

\subsubsec{(1)}~${\mathscr C}=
{{\mathbbm C}}^{\scriptscriptstyle\wedge}$ inside~${{\mathbbm V}}^{({\mathbbm B})}$;

\subsubsec{(2)} there is no nontrivial
${\mathbbm C}$-derivation~in the complex $f$-algebra
${\mathscr C}{\downarrow}$;

\subsubsec{(3)}~each band preserving
endomorphism is a~band projection in~${\mathscr C}{\downarrow}$;

\subsubsec{(4)} there is no band preserving
automorphism other than the identity~in~${\mathscr C}{\downarrow}$.

\subsubsec{(5)}~${\mathbbm B}$ is $\sigma$-distributive.
\endproc

\subsection{5.9}~The above problem was posed by A.~W.~Wickstead
in~\cite{Wic1}. The first example of an unbounded band preserving linear operator
was suggested by Yu.~A. Abramovich,
A.~I.~Veksler, and A.~V.~Koldunov in~\cite{AVK, AVK1}.
Theorem 5.3 combines a~result of
Yu.~A.~Abramovich, A.~I.~Veksler, and A.~V.~Koldunov
\cite[Theorem~2.1]{AVK} and that of P.~T.~N.~McPolin and
A.~W.~Wickstead \cite[Theorem~3.2]{MW}. Theorem 4,7 was obtained
by A.~E.~Gutman \cite{Gut1}; he also found an~example of a~purely
nonatomic locally one-dimensional Dedekind complete vector lattice
(see \cite{Gut6}). Theorem 5.8 belong to A.~G.~Kusraev~\cite{Kus2}.

\newpage
\ssection {6. Order Bounded Operators}

A linear functional  on a~vector space is determined up to a~scalar
from its zero hyperplane. In contrast, a~linear operator is
recovered from its kernel up to a~simple multiplier on a~rather
special  occasion. Fortunately,  Boolean valued analysis
prompts us that some operator  analog of the functional case is valid
for each operator with  target a~Kantorovich space, a~Dedekind
complete vector lattice.
We now proceed along the lines of this rather promising approach.

\subsection{6.1}
Let $E$ be a~vector lattice, and let $F$ be a~$K$-space
with base a~complete Boolean algebra~$\mathbbm B$.
By~4.2, we may assume that
$F$ is a~nonzero space embedded as an order dense ideal in the
universally complete Kantorovich space
$\mathscr R{\downarrow}$ which is the descent of the
reals~$\mathscr R$ inside the separated Boolean valued universe
${\mathbbm V}^{({\mathbbm B})}$  over~$\mathbbm B$.

An~operator $T$~is {\it $F$-discrete\/} if
$[0, T ]=[0, I_F] \circ T$;
i.e., for all
$0 \le S \le T$
there is some $0 \le \alpha \le I_F$
satisfying
$S=\alpha\circ T$.
Let
$L_a^\sim (E, F)$
be the band in~$L^\sim (E, F)$
spanned by
$F$-discrete operators and
$L_d^\sim (E, F)\!:= L _a^\sim (E, F)^\perp$.
By analogy we define
$(E^{\scriptscriptstyle\wedge\sim})_a $
and
$(E^{\scriptscriptstyle\wedge\sim}) _d$.
The members of
$L_d^\sim (E, F)$
are usually called
$F$-{\it diffuse}.

\subsection{6.2}~As usual, we  let $E^{\scriptscriptstyle{\wedge}}$
stand for the standard name of~$E$ in~${\mathbbm V}^{({\mathbbm B})}$.
Clearly, $E^{\scriptscriptstyle{\wedge}}$
is a~vector lattice over~${\mathbbm R}^{\scriptscriptstyle{\wedge}}$
inside~${\mathbbm V}^{({\mathbbm B})}$.
Denote by  $\tau:=T{\uparrow}$ the ascent of~$T$
to~${\mathbbm V}^{({\mathbbm B})}$.
Clearly,  $\tau$ acts from
$E^{\scriptscriptstyle{\wedge}}$ to the ascent~$F{\uparrow}=\mathscr R$
of~$F$ inside the Boolean valued universe ${\mathbbm V}^{({\mathbbm B})}$.
Therefore,
$\tau(x^{\scriptscriptstyle{\wedge}})=Tx$
inside ${\mathbbm V}^{({\mathbbm B})}$ for all $x\in E$,  which means in terms of
truth values that
$[\![\tau:E^{\scriptscriptstyle{\wedge}}\to\mathscr R]\!]=
{\mathbb 1}$ and $(\forall x\in E)\ [\![\tau(x^{\scriptscriptstyle{\wedge}})=Tx]\!]={\mathbb 1}$.

Let
$E^{\scriptscriptstyle\wedge\sim}$
stand for the space of
all order bounded
${\mathbbm R}^{\scriptscriptstyle\wedge}$-linear
functionals from~$E^{\scriptscriptstyle\wedge}$
to~$\mathscr R$.
Clearly,
$E^{\scriptscriptstyle\wedge\sim}\!:=
L^\sim (E^{\scriptscriptstyle\wedge}, \mathscr R)$
is a~$K$-space inside
${\mathbbm V}^{({\mathbbm B})}$.
The descent
$E^{\scriptscriptstyle\wedge\sim}{\downarrow}$
of $E^{\scriptscriptstyle\wedge\sim}$
is a~$K$-space.
Given $S, T \in L^\sim (E, F)$, put $\tau\!:= T{\uparrow}$
and
$\sigma:= S{\uparrow}$.

\subsection{6.3}
\proclaim{Theorem.}
For each $T\in L^{\sim}(E, F)$ the ascent $T{\uparrow}$ of~$T$
is an order bounded ${\mathbbm R}^{\scriptscriptstyle\wedge}$-linear functional
on~$E^{\scriptscriptstyle\wedge}$ inside
${\mathbbm V}^{({\mathbbm B})}$; i.e.,
$[\![T{\uparrow}\in E^{\scriptscriptstyle\wedge\sim }]\!] ={{\mathbb 1}}$.
The mapping $T\mapsto T{\uparrow}$ is a~lattice
isomorphism of  $L^\sim (E, F)$ and $E^{\scriptscriptstyle\wedge\sim}{\downarrow}$.
In particular, the following hold:

{\bf(1)}~$T\ge 0\, \leftrightarrow\,
[\![\,\tau \ge 0\,]\!] ={{\mathbb 1}}$;

{\bf (2)}~$S$ is a~fragment of $T$$\,\leftrightarrow\,
[\![\,\sigma$ is a~fragment of~$\tau\,]\!]={{\mathbb 1}}$;

{\bf (3)}~$T$~is a~lattice homomorphism
if and only if so is $\tau$ inside ${\mathbbm V}^{({\mathbbm B})}$;
%$\leftrightarrow [\![\tau\ \mbox{is a~lattice
%homomorphism}]\!]={\mathbb 1}$;

{\bf (4)}~$T$ is~$F$-diffuse $\,\leftrightarrow\,
[\![\,\tau$ is diffuse $]\!]={{\mathbb 1}}$;

{\bf (5)}~$T \in L_a^\sim (E, F)\,\leftrightarrow\,
[\![\,\tau\in (E ^{\scriptscriptstyle\wedge\sim})_a\,]\!]
={{\mathbb 1}}$;

{\bf (6)}~$T \in L_d^\sim (E, F)\,\leftrightarrow\,[\![\,\tau
\in (E^{\scriptscriptstyle\wedge\sim})_d\,]\!]={{\mathbb 1}}$.
\endproc

Thus, the ascent and descent operations implement a lattice isomorphism of
the vector lattice of all linear order bounded operators from $E$ to $F$ and
the descent of the internal space of all linear order bounded functionals in
the standard name of $E$. This Boolean valued representation does not
preserve order continuity and is not suitable for the study of order
continuous operators. But it may reduce some problems on general order
bounded and positive operators to those on functionals and provide a rather
promising approach.

Consider an instance of this approach. A linear functional on a vector space
is determined up to a scalar from its zero hyperplane. In contrast, a linear
operator is recovered from its kernel up to a simple multiplier on a rather
special occasion. Fortunately, Boolean valued analysis prompts us that some
operator analog of the functional case is valid for each operator with target
a Kantorovich space.

More precisely, since $\tau$, the ascent of an~order bounded operator $T$, is
defined up to a~scalar from $\ker({\tau})$, we infer the following analog of
the Sard Theorem.

\subsection{6.4}
\proclaim{Theorem.}
Let $S$ and $T$ be   linear operators from $E$ to $F$.
Then $\ker(bS)\supset \ker(bT)$
for all $b\in \mathbbm B$   if and only if there is an orthomorphism
$\alpha$ of~$F$ such that   $S=\alpha T$.
\endproc

We see that a~linear operator $T$ is, in a~sense, determined up to
an orthomorphism from the family of the kernels of the
{\it strata\/} $bT$ of~$T$.  This remark opens a~possibility
of studying some properties of $T$ in terms of the kernels of the
strata of~$T$.

\subsection{6.5}
\proclaim{Theorem.}
An order bounded operator $T$
from $E$ to $F$  may be presented as the difference
of some lattice homomorphisms if and only if
the kernel of each stratum $bT$ of~$T$
is a~vector sublattice of~$E$ for all $b\in \mathbbm B$.
\endproc

Straightforward calculations of  truth values show that
$T_+{\uparrow}=\tau_+$
and $T_-{\uparrow}=\tau_-$ inside ${\mathbbm V}^{({\mathbbm B})}$.
Moreover,
$[\![\ker(\tau)$ is a~vector sublattice of~$E^{\scriptscriptstyle{\wedge}}
]\!]= {\mathbb 1}$ whenever so are  $\ker(bT)$ for all $b\in \mathbbm B$.
Since the ascent of a~sum is the sum of
the ascents of the summands, we reduce
the proof of Theorem~6.5 to the case of the functionals
on using 6.3\,(3).

\subsection{6.6}
Recall that a~subspace $H$ of a~vector lattice~$E$ is
a~{\it $G$-space\/} or {\it Grothendieck subspace}  (cp.~\cite{Gro, LinWul})
provided that $H$ enjoys the following property:
$$
(\forall x,y\in H)\ (x\vee y\vee 0 + x\wedge y\wedge 0 \in H).
$$

By simple calculations of truth values we  infer that
$[\![\ker(\tau)$ is a Grothendieck subspace
of~$E^{\scriptscriptstyle{\wedge}}]\!]={\mathbb 1}$
if and only if the kernel of each stratum~$bT$ is
a~Grothendieck subspace of~$E$. We may now assert that the following
appears as  a~result of ``descending''  its scalar analog.

\subsection{6.7}
\proclaim{Theorem.}
The modulus of an order bounded  operator $T: E \to F$
is the sum of some pair of lattice homomorphisms
if and only if the kernel of each stratum
$bT$ of~$T$ with $b\in \mathbbm B$ is a~Grothendieck subspace of
the ambient vector lattice~$E$.
\endproc

To prove the relevant scalar versions of Theorems~6.5 and 6.7, we use
one of the formulas of subdifferential calculus (cp.~\cite{KK3}):

\subsection{6.8}
\proclaim{Decomposition Theorem.}
Assume that $H_1,\dots,H_N$ are cones in a~vector lattice~$E$.
Assume further that $f$ and $g$ are positive functionals on~$E$.
The inequality
$f(h_1\vee\dots\vee h_N)\ge g(h_1\vee\dots\vee h_N)$
holds for all
$h_k\in H_k$ $(k:=1,\dots,N)$
if and only if to each decomposition
of~$g$ into a~sum of~$N$ positive terms
$g=g_1+\dots+g_N$
there is a~decomposition of~$f$ into a~sum of~$N$
positive terms $f=f_1+\dots+f_N$
such that
$f_k(h_k)\ge g_k(h_k)\quad
(h_k\in H_k;\ k:=1,\dots,N)$.
\endproc

\subsection{6.9}
Theorems~6.5 and~6.7 were obtained by S.~S.~Kutateladze
in~\cite{Kut51, Kut52}.
Theorem~6.8 appeared in this form in~\cite{Kut75}.
Note that the sums of lattice homomorphisms were first described
by  S.~J.~Bernau, C.~B.~Huijsmans, and B.~de~Pagter in terms of
$n$-disjoint operators in~\cite{BHP}.
A survey of some conceptually close results on $n$-disjoint operators
is given in~\cite{DOP}.

\ssection{7. Fragments of Positive Operators}

In this section the tools for generating fragments of positive operators and
representation of principal band projections are described. As above, we use
the general method of ascending into a~Boolean-valued universe and descending
the corresponding results for functionals.

\subsection{7.1}~A set of band projections ${\mathscr P}$ in the $K$-space
$L^\sim(E,F)$ {\it generates} the fragments of a positive operator $T\in
L^\sim(E,F)_+$ provided that $Tx^+=\sup\{ (p T) x:\,p\in {\mathscr P}\}$. If the
latter is true for all $T\in L^\sim(E,F)_+$ and $x\in E$ then ${\mathscr P}$ is
said to be a {\it generating set}.  As an~easy example we cite the following.
To each band projection $\pi\in\mathfrak P(E)$ assign the band projection
$\hat{\pi} T\mapsto T\circ\pi$ acting in $L^\sim(E,F)$ and denoted
by ${\mathscr P}^\circ$ the set of all such band projections. It cam be
easily checked that if a~vector lattice $E$ has projection property
then ${\mathscr P}^\circ$ is a~generating set of projections in $L^\sim(E,F)$.

Put ${\mathscr P}^\pi\!:=\{\pi_e:\,e\in E_+\}$ where  $\pi_e$ is defined
as follows:
$$
\gathered
\pi_e Tx=\sup_n T(ne\land x)\quad (x\in E^+,\,T\in L^+(E,F)),\\
\pi_e Tx=\pi_e Tx^+-\pi_e Tx^-\quad(x\in E,\,T\in L^+(E,F)),\\
\pi_eT=\pi_eT^+-\pi_eT^-\quad(T\in L^\sim(E,F)).
\endgathered
$$
Then ${\mathscr P}^\pi$ is a generating set of projections in $L^\sim(E,F)$.

\subsection{7.2}~According to 6.3 the mapping
$T\in L^\sim(E,F)\mapsto T{\uparrow}\in
E^{{\scriptscriptstyle\wedge}\sim}{\downarrow}$
implements an isomorphism between the structures of $L^\sim(E,F)$
and $E^{{\scriptscriptstyle\wedge}\sim}{\downarrow}$.  Therefore,
$T$ is a fragment of $S$ or $T$ is in $S^{\perp\perp}$ if and only
if $T{\uparrow}$ is a fragment of $S{\uparrow}$ or $T{\uparrow}$ is
in $\{S{\uparrow}\}^{\perp\perp}$  inside ${\mathbbm V}^{({\mathbbm B})}$.

The mapping $T{\uparrow}\mapsto (pT){\uparrow}$ $(T\in L^\sim(E,F))$
is extensional for $p\in {\mathscr P}$. By analogy, the ascent
$p{\uparrow}$ is defined to be the band projection in
$E^{{\scriptscriptstyle\wedge}\sim}$ inside
${\mathbbm V}^{({\mathbbm B})}$ acting by the rule
$p{\uparrow} T{\uparrow}=(pT){\uparrow}$ for $T\in L^\sim(E,F)$.

Now, consider the ascent ${\mathscr P}{\uparrow}$ defined as
${\mathscr P}{\uparrow}:=\{p{\uparrow}\mid  p\in {\mathscr P} \}{\uparrow}$.
Obviously, ${\mathscr P}$ generates the fragments of $T$ if and only if
${\mathscr P}{\uparrow}$ generates the fragments of $T{\uparrow}$
inside ${\mathbbm V}^{({\mathbbm B})}$.

Given a set $A$ in a $K$-space, we denote by $A^\vee$ the union of
$A$ and the suprema of all nonempty finite subsets of~$A$. The symbol
$A^{({\uparrow})}$ denotes the result of adjoining to $A$ the
suprema of all increasing nonempty nets in $A$. The symbols
$A^{({\uparrow}{\downarrow})}$ and
$A^{({\uparrow}{\downarrow}{\uparrow})}$
are interpreted in a natural way.

\subsection{7.3}~Let ${\mathbbm P}$ be a dense subfield of
${\mathbbm R}$ and $E$ be a~vector lattice over ${\mathbbm P}$. Denote
by $E^\sim\!:=L^\sim(E,{\mathbbm R})$ the vector lattice of ${\mathbbm
P}$-linear functionals in $E$. Fix some set ${\mathscr P}$ of band projections
and  the corresponding set ${\mathscr P}(f)\!:=\{ pf:\, p\in{\mathscr P}\}$
of the fragments of a positive functional $f\in E^\sim$.

\Proclaim{Theorem.}
For positive functionals $f,g\in E^\sim$ the following  are true:

\subsubsec{(1)}~${\mathscr P}$ generates the fragments of $f$ if and only if
${\mathscr P} (f)^{\vee({\uparrow}{\downarrow}{\uparrow})}=\mathfrak E(f)$;

\subsubsec{(2)}~if ${\mathscr P}$ is generating then $g\in \{f\}^{\perp\perp}$
if and only if for any $x\in E_+$ and $0<\varepsilon\in {\mathbbm R}$
there exists $0<\delta\in{\mathbbm R}$ such that  $pf(x)\le \delta$
implies $pg(x)\le\varepsilon$ for every $p\in {\mathscr P}$;

\subsubsec{(3)}~if ${\mathscr P}$ is generating then for the principal band
projection $\pi_f$ onto~$\{f\}^{\perp\perp}$ the representations hold:
$$
\pi_f g(x)=\sup_{\varepsilon>0}\inf\{pg(x):\
p^\perp f(x)\le \varepsilon,\ p\in {\mathscr P} \}.
$$
\endproc

\subsection{7.4}~
\proclaim{Theorem.}
A set ${\mathscr P}$ of band projections in $L^\sim(E,F)$ generates
the fragments of $T\in L^\sim(E,F)_+$ if and only if
${\mathscr P} (T)^{\vee({\uparrow}{\downarrow}{\uparrow})}=\mathfrak E(T)$.
\endproc

\subsection{7.5}~
\proclaim{Theorem.}
If ${\mathscr P}$ is a generating set of projections in $L^\sim(E,F)$
then for positive operators $S,T\in L^\sim(E,F)$ the relation
$T\in \{S\}^{\perp\perp}$ holds if and only if for every
$e\in E_+$ and $0<\varepsilon\in{\mathbbm R}$ there exits
$0<\delta\in F$, $[Se]\le [\delta]$, $\delta\le Se$, such that
$\pi pSe\le\delta$ implies $\pi pTe\le\varepsilon Te$ for all
$\pi\in\mathfrak P(F)$ and $p\in{\mathscr P}$.
\endproc

\subsection{7.6}
\proclaim{Theorem.}
Let $E$ be an~arbitrary vector lattice, $F$ be a $K$-space, and
${\mathscr P}$ be a  generating set of projections. Then for the band
projection $T_S$ of $T$ onto $\{S\}^{\perp\perp}$ the
representations are valid:
$$
\gathered
(T-T_S)e=\inf_{0<\varepsilon \in\Bbb R}\sup\{\pi pTe: \pi\in\mathfrak P(F),\,
p\in{\mathscr P},\,\pi pSe\le\varepsilon Se\},\\
(T_S)e=
\sup_{0<\varepsilon \in\Bbb R}\inf\{(\pi p)^\perp Te: \pi\in\mathfrak P(F),\,
p\in{\mathscr P},\,\pi pSe\le\varepsilon Se\}.
\endgathered
$$
\endproc

\subsection{7.7}~The concept of a~generating set of projections
as well as Theorem 7.4 belongs to S.~S.~Kutateladze
\cite{Kut}. In 7.4 every fragment of a~positive operator is
obtained from its simpler fragments by up and down procedures.
Similar assertions are often referred to as {\it up-down
theorems}. The first up-down theorem was established by
B.~de~Pagter \cite{Pag3} (also see \cite{AB1}, \cite{AB}).
However, it involved two essential constraints: $F$ should admit
a~total set of \hbox{$o$}-continuous functionals, and $E$ must
be order complete (or at least possess the principal projection
property). The first constraint was eliminated in~\cite{KS}
and the second, in~\cite{AKK}. Of course, a few up-down
theorems can be deduced from 7.4 by specifying generating sets
(see  \cite{DOP} for details). Theorems 7.5 and 7.6 are improved
versions of the corresponding results of~\cite{Kut}.

\ssection{8. Boolean Valued Banach Spaces}

In this section we discuss the~transfer principle of Boolean valued analysis
in regard to lattice-normed spaces. It turns out that
the~interpretation of a~Banach space inside
an~arbitrary Boolean valued model is a~Banach--Kantorovich space.
Conversely, the~universal completion of each lattice-normed space
becomes a~Banach space on ascending in a~suitable Boolean valued model.
This open up an opportunity to transfer the available theorems on Banach spaces to
analogous results on lattice-normed spaces by the technique of Boolean valued analysis.

\subsection{8.1}
Consider a~vector space $X$ and a~real vector
lattice $E$. Note that all vector lattices under consideration are
assumed Archimedean.
An~{\it $E$-valued norm\/}  is a~mapping
$\[{\cdot}\]:X\to E_+$ such that

\subsubsec{(1)}~$\[x\]=0 \Longleftrightarrow x=0  \quad (x\in X)$;

\subsubsec{(2)}~$\[\lambda x\]=|\lambda|\[x\] \quad (\lambda\in{\mathbbm R}$, $x\in X)$;

\subsubsec{(3)}~$\[x+y\] \leq \[x\]+\[y\] \quad (x,y\in X)$.

A~vector norm is  {\it decomposable\/} if

\subsubsec{(4)}~for all $e_1, e_2 \in E_+$ and $x \in X$, from
$\[x\]=e_1+e_2$ it follows that there exist $x_1,x_2 \in X$ such that
$x=x_1+x_2$ and $\[x_k\]=e_k$ $(k \!:= 1,2)$.

If~(4) is valid only for disjoint
$e_1, e_2 \in E_+$ then  the~norm is  {\it $d$-decomposable}.
A~triple $\bigl( X,\[{\cdot}\],E\bigr)$ as well as
briefer versions is a~{\it lattice-normed space\/}
over $E$ whenever $\[{\cdot}\]$ is an~$E$-valued norm on~$X$.

\subsection{8.2}
By a~{\it Boolean algebra of projections\/}
in a~vector space $X$ we mean a~set ${\mathscr B}$ of commuting
idempotent linear operators in~$X$.
Moreover, the~Boolean operations have the~following form:
$$
\gathered
\pi \wedge \rho\!:=\pi \circ \rho=\rho \circ \pi, \quad
\pi \vee \rho=\pi + \rho - \pi \circ \rho,
\\
\pi^\ast=I_x-\pi \quad (\pi,\rho \in {\mathscr B}),
\endgathered
$$
and the~zero and identity operators in $X$ serve as
the~zero and unity of the Boolean algebra ${\mathscr B}$.

\Proclaim{}
Suppose that $E$ is a~vector lattice with the projection
property and
$E=\[X\]^{\perp\perp}\!:=\{\[x\]:\,x\in X\}^{\perp\perp}$.
If $(X,E)$ is a~$d$-decomposable lattice-normed space
then there exists a~complete Boolean algebra ${\mathscr B}$ of band
projections in $X$ and an~isomorphism $h$ from ${\mathfrak P}(E)$
onto ${\mathscr B}$ such that
$$
b\[x\]=\[h(b)x\]\quad\bigl(b\in{\mathfrak P}(E),\ x\in X\bigr).
$$
\endproc

We identify the Boolean algebras ${\mathfrak P}(E)$ and
${\mathscr B}$ and write
$\pi \[x\]=\[\pi x\]$
for all $x\in X$ and $\pi\in{\mathfrak P}(E)$.

\subsection{8.3}
A~net $(x_\alpha)_{\alpha\in {\roman A}}$
in $X$ is {\it $bo$-convergent\/} to $x\in X$ (in symbols:
$x=\bolim x_\alpha$) if  $(\[x-x_\alpha\])_{\alpha\in {\roman A}}$
is $o$-convergent to zero. A~lattice-normed space $X$ is
{\it $bo$-complete} if each net $(x_\alpha)_{\alpha\in{\roman A}}$
is $bo$-convergent to  some element of~$X$ provided that
$(\[x_\alpha -x_\beta\])_{(\alpha,\beta)\in {\roman A}\times {\roman A}}$
is $o$-convergent to zero.
A decomposable $bo$-complete lattice-normed space $(X,\[\cdot\], E)$
is called a~{\it Banach--Kantoro\-vich space}. If  $E$
is a~universally complete Kantorovich space then $X$ is also referred to as
universally complete.
By a~{\it universal completion\/} of a~lattice-normed space $(X,E)$
we mean a~universally complete Ba\-nach--Kantorovich space $(Y,m(E))$
together with a~linear isometry $\imath :X \rightarrow Y$ such that
each universally complete $bo$-complete subspace of $(Y,m(E))$
containing $\imath (X)$ coincides with $Y$. Here $m(E)$ is a~universal
completion~of~$E$.

\subsection{8.4}
\proclaim{Theorem.}
Let $({\mathscr X},\,\|\cdot\|)$ be a~Banach space inside~${\mathbbm V}^{({\mathbbm B})}$.
Put $X\!:={\mathscr X}{\downarrow}$ and
$\[{\cdot}\]:=\|\cdot\|{\downarrow}({\cdot})$.  Then
$\bigl(X,\[{\cdot}\],{\mathscr R}{\downarrow}\bigr)$ is a~universally
complete Banach--Kantorovich space. Moreover, $X$ can be endowed
with the~structure of a~faithful unitary module over the~ring
$\Lambda:={\mathscr C}{\downarrow}$ so that $\[ax\]=|a|\,\[x\]$ and
$b\leq [\![\,x=0\,]\!]\,\leftrightarrow\,\chi(b)x=0$
for all $a\in {\mathscr C}{\downarrow}$, $x\in X$, and $b\in {\mathbbm B}$,
where $\chi$ is an~isomorphism of\/ ${\mathbbm B}$ onto ${\mathfrak P}(X)$.
\endproc

\subsection{8.5}
\proclaim{Theorem.}
To each lattice-normed space $(X,\[\cdot\])$, there exists a
unique  Banach space (up to a~linear isometry) ${\mathscr X}$
inside ${\mathbbm V}^{({\mathbbm B})}$, with ${\mathbbm B}\simeq
{\mathfrak B}\,\bigl(\[X\]^{\perp\perp}\bigr)$, such that the~descent
${\mathscr X}{\downarrow}$ of~$\mathscr X$ is a~universal completion of~$X$.
\endproc

As~in~4.1, we call  $x\in X$ and $y\in Y$  {\it disjoint\/}
and write $x\perp y$ whenever $\[x\]\land\[y\]=\nobreak0$. Let $X$
and $Y$ be Banach--Kantorovich spaces over some $K$-space~$G$.

An~operator~$T$ is  {\it band preserving\/} if
$x\perp y$ implies $Tx\perp y$ for all $x\in X$ and $y\in Y$.
Denote by ${\mathscr L}_G (X,Y)$ the space of all band preserving
operators $T:X\to Y$ that send all norm-$o$-bounded sets into
norm-$o$-bounded sets.

\subsection{8.6}
\proclaim{Theorem.}
Let ${\mathscr X}$ and ${\mathscr Y}$ be Boolean valued representations for
Banach--Kantoro\-vich spaces $X$ and $Y$ normed by some universally
complete $K$-space~$G\!:={\mathscr R}{\downarrow}$. Let
${\mathscr L}^{\mathbbm B}({\mathscr X},{\mathscr Y})$ be the space of bounded linear
operators from ${\mathscr X}$ into ${\mathscr Y}$ inside ${\mathbbm V}^{({\mathbbm B})}$,
where ${\mathbbm B}\!:={\mathfrak B}(G)$. The descent and ascent mappings
(for operators) implement linear isometries between the
lattice-normed spaces ${\mathscr L}_G (X,Y)$ and
${\mathscr L}^{\mathbbm B}({\mathscr X},{\mathscr Y}){\downarrow}$.
\endproc

\subsection{8.7}~The concept of lattice-normed space was
suggested  by L.~V.~Kantorovich in~1936
\cite{Ka1}. It is worth stressing that \cite{Ka1}
is the fist article with the~unusual decomposability
axiom for an~abstract norm. Paradoxically, this axiom  was often omitted
as inessential in the further papers by other authors.
The profound importance of 8.1\,(4) was revealed by
Boolean valued analysis. The connection between the decomposability and
existence of a~Boolean algebra of projections in a~lattice-normed space
was discovered in~\cite{K1, KVD}. The theory
of lattice-normed spaces and dominated operators is
set forth in~\cite{DOP}. As regards the Boolean valued approach,
see \cite{IBA}.

\ssection{9. Boolean Valued Order Continuous Functionals}

We now address the class of $o$-continuous order bounded
operators that turn into $o$-continuous functionals on ascending to
a suitable Boolean valued model.

\subsection{9.1}
Assume that a~lattice-normed
space $X$ is simultaneously a~vector lattice.
The norm $\[\cdot\]:X\to E_+$  of $X$ is
{\it monotone\/} if from
$|x|\le|y|$ it follows that $\[x\]\le \[y\]$ ($x,y\in X$).
In this event,  $X$~is a~{\it lattice-normed vector lattice}.
Moreover, if $X$~is a~Banach--Kantorovich space then
$X$ is called a~{\it Banach--Kantorovich lattice}.

We say that the norm $\[\cdot\]$ in $X$ is {\it additive }
if $\[x+y\]=\[x\]+\[y\]$ for all  $x,y\in X_+$; it is
{\it order semicontinuous\/}
or $o$-{\it semicontinuous\/} for short if
$\sup\[x_\alpha\]=\[\sup x_\alpha \]$ for each increasing net
$(x_\alpha)\subset X$ with the least upper bound $x\in X$; and
it is {\it order continuous\/} or $o$-{\it continuous\/} if
$\inf\[x_\alpha\]=0$ for every decreasing net
$(x_\alpha)\subset X$ with $\inf_\alpha x_\alpha=0$.

The Boolean valued interpretation of Banach--Kantorovich lattices
proceeds along the lines of the previous section.

\subsection{9.2}
\proclaim{Theorem.}
Let $(X,\[\cdot\])$ be a~Banach--Kantorovich space
and let $({\mathscr X},\|\cdot\|)\in{{\mathbbm V}}^{({\mathbbm B})}$ stand for
its Boolean valued realization. Then

\subsubsec{(1)}~$X$ is a~Banach--Kantorovich lattice
if and only if
${\mathscr X}$~is a~Banach lattice
inside ${{\mathbbm V}}^{({\mathbbm B})}$;

\subsubsec{(2)}~$X$
is an~order complete Banach--Kantorovich lattice
if and only if
${\mathscr X}$ is~an order complete Banach lattice
inside ${{\mathbbm V}}^{({\mathbbm B})}$;

\subsubsec{(3)}~the norm $\[\cdot\]$ is $o$-continuous
(order semicontinuous, monotone complete, or additive)
if and only if the norm $\|\cdot\|$ is $o$-continuous
(order semicontinuous, monotone complete, or
additive) inside ${{\mathbbm V}}^{({\mathbbm B})}$.
\endproc

\subsection{9.3}
Let $E$ be a~vector lattice, let $F$ be some
$K$-space, and let $T$ be a~positive operator from $E$ to~$F$.

Say that $T$ possesses the {\it Maharam property} if, for
all $x\in E_+$ and $0\le f\le Tx\in F_+$, there  is some
$0\le e\le x$ satisfying $f=Te$. An $o$-continuous
positive operator with the Maharam property is  a~{\it Maharam
operator}.

Observe that $T\in L(E, F)_+$ possesses the Maharam
property if only if the equality $T\,(\,[0,x]\,)=[0,Tx]$
holds for all $x\in E_+$. Thus, a~Maharam operator is exactly
an~$o$-continuous order-interval preserving positive operator.

Let $T$ be an essentially positive operator from
$E$ to~$F$ enjoying the Maharam property. Put
$\[e\]:=T(|x|)$ $(e\in E)$. Then $(E,\[\cdot\])$ is a~disjointly
decomposable lattice-normed space over~$F$.

Put $F_T\!:=\{T(|x|):\,x\in E\}^{\perp\perp}$, and let
${\mathscr D}_m\,(T)$ stand for the greatest order dense ideal
of the universal completion $m(E)$ of~$E$ among those to which
$T$ can be extended by $o$-continuity. In other words,
$z\in{\mathscr D}_m\,(T)$ if and only if
$z\in m(E)$ and the set $\{T(x):\,x\in E,\,0\leq x\leq|z|\}$
is bounded in~$F$. In this event there exists a~minimal extension
of~$T$ to~${\mathscr D}_m(T)$ presenting an $o$-continuous positive
operator.

\Proclaim{}
Let $E$  and $F$ be some $K$-spaces, and
let $T:E\to F$ be a~Maharam operator. Put
$X\!:={\mathscr D}_m(T)$ and $\[x\]\!:=\Phi(|x|)$ $(x\in X)$,
where $\Phi$ is an $o$-continuous extension of~$T$ to~$X$.
Then $(X,\[\cdot\])$ is a~Banach--Kantorovich lattice
whose norm is $o$-continuous and additive.
\endproc

\subsection{9.4}
\proclaim{Theorem.}
Let $X$ be an~arbitrary $K$-space and let $E$ be a~universally
complete $K$-space ${\mathscr R}{\downarrow}$. Assume that
$\Phi:X\rightarrow E$ is a~Maharam operator such that
$X=X_\Phi={\mathscr D}_m\,(\Phi)$ and $E=E_\Phi$. Then there are
elements ${\mathscr X}$ and $\varphi$ in~${\mathbbm V}^{({\mathbbm B})}$ satisfying

\subsubsec{(1)}~$[\![{\mathscr X}$ is a~$K$-space, $\varphi:
{\mathscr X}\rightarrow{\mathscr R}$ is a~positive $o$-continuous functional, and
${\mathscr X}={\mathscr X}_\varphi={\mathscr D}_m\,(\varphi)]\!]={\mathbb 1};$

\subsubsec{(2)}~if $X'\!:={\mathscr X}{\downarrow}$ and
$\Phi'=\varphi{\downarrow}$ then $X'$ is a~$K$-space and
$\Phi':X'\rightarrow E$ is a~Maharam operator;

\subsubsec{(3)}~there is a~linear and lattice isomorphism $h$
from $X$ onto $X'$ such that $\Phi=\Phi'\circ h$;

\subsubsec{(4)}~for a~linear operator $\Psi$, the containment
$\Psi\in\{\Phi\}^{\perp\perp} $ is true if and only if there is
$\psi\in{\mathbbm V}^{({\mathbbm B})}$ such that  $\psi\in\{\varphi\}^{\perp\perp}$
inside ${\mathbbm V}^{({\mathbbm B})}$ and $\Psi=(\psi{\downarrow})\circ h$.
\endproc

Theorem 9.4 enables us to claim that
each fact about $o$-continuous positive linear functionals
in $K$-spaces
has a~parallel version for Maharam operators which can be
revealed by using 9.4. For instance, we state the abstract

\subsection{9.5}
\proclaim{Radon--Nikod\'ym Theorem.}
Let $E$ and $F$ be $K$-spaces. Assume further that $S$
and $T$ are $o$-continuous positive operators from $E$ to $F$,
with~$T$ enjoying the Maharam property. Then the following
are equivalent:

\subsubsec{(1)}~$S\in \{T\}^{\perp\perp}$;

\subsubsec{(2)}~$Sx\in\{Tx\}^{\perp\perp}$ for all $x\in E_+$;

\subsubsec{(3)}~there is an extended orthomorphism
$0\le\rho\in\Orth^\infty(E)$ satisfying
$Sx=T(\rho x)$
for all
$x\in E$ such that $\rho x\in E$;
%{\mathscr D}(\rho )$;

\subsubsec{(4)}~there is a~sequence of orthomorphisms
$(\rho_n)\subset\Orth(E)$
such that
$Sx=\sup_n T(\rho_nx)$
for all $x\in E$.
\endproc

\subsection{9.6}
A~brief description for  Maharam's
approach to studying positive operators in the~spaces
of measurable functions and the main results in this area
are collected in~\cite{Mah5}. W.~A.~J.~Luxemburg and A.~R.~Schep~\cite{LS}
extended a~portion of Maharam's theory on
the~Radon--Nikod\'ym Theorem to the~case of positive operators in
vector lattices.

Theorem~9.2 and 9.4 were obtained by A.~G.~Kusraev
\cite{K12}
and Theorem~9.5, by W.~A.~J. Luxemburg and A.~R.~Schep~\cite{LS}.
About various applications of the above results on Maharam operators
and some extension of this theory to sublinear and convex
operators see \cite{KVD, DOP,  KK3,  IBA}.

\ssection{10. Spaces with Mixed Norm}

The definitions of various objects of functional analysis rest often
on some blending of  the norm and order properties. Among these are
listed the spaces with mixed
norm and the classes of linear operators between them.

\subsection{10.1}
If~$(X,E)$ is a~lattice-normed space whose
{\it norm lattice\/} $E$ is a~Banach lattice. Since, by definition, $\[x\]\in E$ for
$x\in X$, we may introduce the~{\it mixed norm\/} on~$X$
by the~formula
$$
\tvert x\tvert :=\|\[x\]\|\quad (x\in X).
$$
In~this situation, the~normed space $(X,\tvert \cdot\tvert )$
is called a~{\it space with mixed norm}. A~{\it Banach space
with mixed norm\/} is a~pair $(X,E)$ with~$E$  a~Banach
lattice and~$X$  a~$br$-complete lattice-normed space with
$E$-valued norm. The~following proposition justifies this
definition.

\Proclaim{}
Let~$E$ be a~Banach lattice. Then $(X,\tvert \cdot\tvert )$ is
a~Banach space if and only if the lattice-normed space~$(X,E)$
is  relatively uniformly complete.
\endproc

\subsection{10.2}
Let $\Lambda$ be the {\it bounded part\/} of the
universally complete $K$-space ${\mathscr R}{\downarrow}$, i.e.
$\Lambda$ is the~order-dense ideal in ${\mathscr R}{\downarrow}$
generated by the order unity
${{\mathbb 1}}\!:=1^{\scriptscriptstyle\wedge}\in{\mathscr R}{\downarrow}$.
Take a~Banach space ${\mathscr X}$ inside ${\mathbbm V}^{({\mathbbm B})}$. Put
$$
\gathered
{\mathscr X}{\Downarrow}:=\{x\in{\mathscr X}{\downarrow}:\,\[x\]\in\Lambda \},
\\
\tvert x\tvert \!:=\|\[x\]\|_\infty:=
\inf\{0<\lambda\in{\mathbbm R}:\,\[x\]\leq\lambda{\mathbb 1}\}.
\endgathered
$$
Then ${\mathscr X}{\Downarrow}$ is a~Banach--Kantorovich space called
the {\it bounded descent\/} of ${\mathscr X}$. Since $\Lambda$ is an~order
complete $AM$-space with unity, ${\mathscr X}{\Downarrow}$
is a~Banach space with mixed norm over~$\Lambda$.

Thus, we~came to the~following natural question: Which Banach
spaces are linearly isometric to the bounded descents of internal
Banach spaces? The answer is given in terms of ${\mathbbm B}$-cyclic
Banach spaces.

\subsection{10.3}
Let $X$ be a~normed space. Suppose that
${\mathscr L}(X)$ has a~complete Boolean algebra of norm one projections
${\mathscr B}$ which is isomorphic to ${\mathbbm B}$. In this event we will identify
the Boolean algebras ${\mathscr B}$ and ${\mathbbm B}$, writing ${\mathbbm B}\subset{\mathscr L}(X)$.
Say that $X$ is a~{\it normed ${\mathbbm B}$-space\/} if ${\mathbbm B}\subset{\mathscr L}(X)$
and for every partition of unity $(b_\xi)_{\xi\in\Xi}$ in ${\mathbbm B}$ the
two conditions hold:

\subsubsec{(1)}~If $b_\xi x=0\ (\xi\in\Xi)$ for some
$x\in X$ then $x=0$;

\subsubsec{(2)}~If $b_\xi x=b_\xi x_\xi\ (\xi\in\Xi)$
for $x\in X$ and a~family $(x_\xi)_{\xi\in\Xi}$ in $X$
then
$$
\|x\|\leq\sup\{\|b_\xi x_\xi\|:\xi\in\Xi\}.
$$

Given a~partition of unity $(b_\xi)$, we refer to $x\in X$
satisfying the condition $(\forall\,\xi\in\Xi)\,b_\xi x=b_\xi x_\xi$
as a~{\it mixing\/} of $(x_\xi)$ by $(b_\xi)$. If (1) holds
then there is a~unique mixing $x$ of $(x_\xi)$ by $(b_\xi)$.
In these circumstances we  naturally call $x$ the
{\it mixing\/} of $(x_\xi)$ by $(b_\xi)$. Condition~(2) maybe
paraphrased as follows: The unit ball $U_X$ of $X$ is closed under
mixing.

A normed ${\mathbbm B}$-space $X$ is ${\mathbbm B}$-{\it cyclic\/} if we may find in
$X$ a~mixing of each norm-bounded family by each partition of unity
in~${\mathbbm B}$. It is easy to verify that $X$ is a~${\mathbbm B}$-cyclic normed space
if and only if, given a~partition of unity $(b_\xi)\subset {\mathbbm B}$ and
a family $(x_\xi)\subset U_X$, we may find a~unique element
$x\in U_X$ such that $b_\xi x=b_\xi x_\xi$ for all $\xi$.

A~linear operator (linear isometry) $S$ between normed
${\mathbbm B}$-spaces is {\it ${\mathbbm B}$-linear\/}
({\it ${\mathbbm B}$-isometry\/}) if $S$ commutes with the
projections in ${\mathbbm B}$; i.e.,
$\pi\circ S=S\circ\pi$ for all $\pi\in{\mathbbm B}$.
Denote by ${\mathscr L}_{{\mathbbm B}}(X,Y)$ the set of all bounded
${\mathbbm B}$-linear operators from $X$ to $Y$. We call
$X^{\scriptscriptstyle\#}\!:={\mathscr L}_{{\mathbbm B}}(X,{\mathbbm B}({\mathbbm R}))$
the ${\mathbbm B}$-{\it dual\/} of $X$. If $X^{\scriptscriptstyle\#}$ and
$Y$ are  ${\mathbbm B}$-isometric to each other then we say that
$Y$ is a~{\it ${\mathbbm B}$-dual space\/} and $X$
is a~{\it ${\mathbbm B}$-predual\/}
of $Y$.

\subsection{10.4}
\proclaim{Theorem.}
A Banach space $X$ is linearly isometric to the bounded descent
of some Banach space ${\mathscr X}$ inside ${\mathbbm V}^{({\mathbbm B})}$ (called
a~{\it Boolean valued representation of\/} $X$) if and only if $X$
is ${\mathbbm B}$-cyclic.
If $X$ and $Y$ are ${\mathbbm B}$-cyclic Banach spaces and
${\mathscr X}$ and ${\mathscr Y}$ stand for some Boolean valued
representations of $X$ and $Y$, then the space
${\mathscr L}_{{\mathbbm B}}(X,Y)$ is ${\mathbbm B}$-isometric to the bounded
descent of the internal space ${\mathscr L}({\mathscr X},{\mathscr Y})$ of all
bounded linear operators from ${\mathscr X}$ to ${\mathscr Y}$.
\endproc

\subsection{10.5}
Let $\Lambda$ be a~Stone algebra with unity
${\mathbb 1}$ ($=$~an~order complete complex $AM$-space with strong
order unity ${\mathbb 1}$ and uniquely defined multiplicative structure)
and consider a~unitary $\Lambda$-module~$X$. The mapping
$\langle\cdot\,|\,\cdot\rangle:X\times X\rightarrow\Lambda$
is a~{\it $\Lambda$-valued inner product}, if for all $x,  y,  z\in X$
and $a\in\Lambda$ the following are satisfied:
\medskip

\subsubsec{(1)}~$\langle x\,|\,x\rangle\geq{\mathbb 0};\
\langle x\,|\,x\rangle={\mathbb 0}\leftrightarrow x={\mathbb 0};$
\medskip

\subsubsec{(2)}~$\langle x\,|\,y\rangle =\langle y\,|\,x\rangle^\ast;$
\medskip

\subsubsec{(3)}~$\langle ax\,|\,y\rangle =a\langle x\,|\,y\rangle;$
\medskip

\subsubsec{(4)}~$\langle x+y\,|\,z\rangle=\langle x\,|\,z\rangle +
\langle y\,|\,z\rangle$.
\medskip

Using a~$\Lambda$-valued inner product, we may introduce the norm
of $x\in X$ by
$$
\tvert x\tvert\!:=\sqrt{\|\langle x|x\rangle\|}
$$
and the decomposable vector norm of $x\in X$ by
$$
\[x\]\!:=\sqrt{\langle x|x\rangle}.
$$
Obviously, $\tvert x\tvert =\bigl\|\[x\]\bigr\|$
for all $x\in X$,
and so $X$ is a~space with~mixed norm.

\subsection{10.6}
Let $X$ be a~$\Lambda$-module with an~inner
product $\langle\cdot\,|\,\cdot\rangle:X\times X\rightarrow\Lambda$.
If $X$ is complete with respect to the mixed norm $\tvert\cdot\tvert$
then $X$ is called a~$C^\ast$-{\it module\/} over $\Lambda$. It can be proved
(see \cite{DOP}) that  for a~$C^\ast$-module $X$ the pair
$(X,\tvert\cdot\tvert)$  is a~${\mathbbm B}$-cyclic Banach space if and only if
$(X,\[\cdot\])$ is a~Banach--Kantorovich space over $\Lambda$.
If a~unitary $C^\ast$-module satisfies one of these equivalent
conditions then it is called a~{\it Kaplansky--Hilbert module.}

\subsection{10.7}
\proclaim{Theorem.}
The bounded descent of a~Hilbert space  in
${\mathbbm V}^{({\mathbbm B})}$ is a~Kaplan\-sky--Hilbert module over
the Stone algebra ${\mathscr C}{\Downarrow}$. Conversely, if
$X$ is a~Kaplansky--Hilbert module over ${\mathscr C}{\Downarrow}$,
then there is a~Hilbert space ${\mathscr X}$ in ${\mathbbm V}^{({\mathbbm B})}$
whose bounded descent is unitarily equivalent with $X$.
This space is unique up to  unitary equivalence inside~${\mathbbm V}^{({\mathbbm B})}$.
\endproc

\subsection{10.8}
\proclaim{Theorem.}
Let ${\mathscr X}$ and ${\mathscr Y}$ be Hilbert
spaces inside ${\mathbbm V}^{({\mathbbm B})}$.
Suppose that $X$ and $Y$ are the
bounded descents of ${\mathscr X}$ and ${\mathscr Y}$.
Then the space ${\mathscr L}_{{\mathbbm B}}(X,Y)$ of all
${\mathbbm B}$-linear
bounded operators is a~${\mathbbm B}$-cyclic Banach space
${\mathbbm B}$-isometric to the bounded descent
of the internal Banach space ${\mathscr L}^ {\mathbbm B}({\mathscr X},{\mathscr Y})$ of bounded linear
operators from ${\mathscr X}$ to ${\mathscr Y}$.
\endproc

\subsection{10.9}
Boolean valued analysis approach gives rise
to an~interesting concept of cyclically compact operator in
a~Banach ${\mathbbm B}$-space \cite[8.5.5]{DOP}. Without plunging into
details we formulate a~result on the general form of cyclically
compact operators in Kaplansky--Hilbert modules.

\Proclaim{Theorem.}
Let $X$ and $Y$ be Kaplansky--Hilbert modules over a~Stone
algebra $\Lambda$ and let $T$ be a~cyclically compact operator
from $X$ to $Y$. There are orthonormal families $(e_k)_{k\in{\mathbbm N}}$
in~$X$, $(f_k)_{k\in{\mathbbm N}}$ in~$Y$, and a~family
$(\mu_k)_{k\in{\mathbbm N}}$ in~$\Lambda $ such that the following hold:

\subsubsec{(1)}
$\mu_{k+1}\le\mu_k$
$(k\in{\mathbbm N})$
and
$\olim_{k\to\infty}\mu_k=0;$

\subsubsec{(2)}
there exists a~projection
$\pi_\infty$
in
$\Lambda$
such that
$\pi_\infty\mu_k$
is a~weak order unity in
$\pi_\infty\Lambda$
for all
$k\in{\mathbbm N};$

\subsubsec{(3)}
there exists a~partition
$(\pi_k)_{k=0}^\infty$
of the projection
$\pi_\infty^\perp$
such that
$\pi_0\mu_1=0$,
$\pi_k\le\mu_k$,
and
$\pi_k\mu_{k+1}=0$ for all
$k\in {\mathbbm N};$

\subsubsec{(4)}~the representation is valid
$$
\gathered
T=\pi_\infty\,\bosum\limits_{k=1}^\infty\mu_k e^{\scriptscriptstyle\#}_k
\otimes f_k
\\
+ \bosum\limits_{n=1}^\infty\pi_n\sum_{k=1}^n\mu _k
e^{\scriptscriptstyle\#}_k\otimes f_k.
\endgathered
$$
\endproc

\subsection{10.10}
The  bounded descent of~10.2 appeared in
the research by G.~Takeu\-ti into von Neumann algebras and
$C^\ast$-algebras within Boolean valued models \cite{Tak3, Tak4} and
in the research by M.~Ozawa into Boolean valued
interpretation of the theory of Hilbert spaces~\cite{Ozawa1}.
Theorems~10.4 and 10.9 were obtained by A.~G.~Kusraev in
\cite{K1, KVD, DOP}. Theorems 10.7 and~10.8
were proved by M.~Ozawa \cite{Ozawa1}.

%\newpage
\ssection{11. Banach Algebras}

The possibility of applying Boolean valued analysis to operator
algebras rests on the following observation:  If the center of
an algebra is properly qualified and perfectly located then it
becomes a~one-dimensional subalgebra after ascending in a
suitable Boolean valued universe. This might lead to a~simpler
algebra. On the other hand, the transfer principle implies that
the scope of the formal theory of the initial algebra is the
same as that of its Boolean valued representation.

\subsection{11.1}
An~{\it $AW^\ast$-algebra\/}
is a~$C^\ast$-algebra presenting a~Baer $\ast$-algebra.
More explicitly, an~$AW^\ast$-algebra is a~$C^\ast$-algebra
$A$ whose every {\it right annihilator\/}
$M^\perp\!:=\{y\in A:\,{(\forall x\in M)}\ xy=0\}$
has the form~$pA$, with $p$ a~projection. A~{\it projection\/}~$p$
is a~hermitian ($p^*=p$) idempotent ($p^2=e$) element.
An element $z\in A$ is said to be {\it central\/}
if it commutes with every member of~$A$. The
{\it center\/} of an $AW^*$-algebra $A$ is the set
${\mathscr Z}(A)$ of all central elements. Clearly,
${\mathscr Z}(A)$ is a~commutative $AW^*$-subalgebra of~$A$, with
$\lambda{\mathbb 1}\in {\mathscr Z}(A)$ for all $\lambda\in{\mathbbm C}$.
If ${\mathscr Z}(A)=\{\lambda{\mathbb 1}:\,\lambda\in{\mathbbm C}\}$
then the $AW^*$-algebra $A$ is called an~$AW^*$-{\it factor}.

The symbol ${\mathfrak P}(A)$ stands for the set of all projections
of an involutive algebra $A$. Denote the set of all central projections
by ${\mathfrak P}_c(A)$.

\subsection{11.2}
\proclaim{Theorem.}
Assume that ${\mathscr A}$ is an~$AW^\ast$-algebra inside ${\mathbbm V}^{({\mathbbm B})}$ and $A$~is the bounded descent of~${\mathscr A}$. Then $A$
is also an $AW^\ast$-algebra and, moreover, ${\mathfrak P}_c(A)$  has
an order-closed subalgebra isomorphic with ${\mathbbm B}$.  Conversely, let
$A$ be an~$AW^\ast$-algebra such that ${\mathbbm B}$ is an order-closed
subalgebra of the Boolean algebra ${\mathfrak P}_c(A)$. Then there is
an~$AW^\ast$-algebra ${\mathscr A}$ in ${\mathbbm V}^{({\mathbbm B})}$ whose bounded
descent is $\ast$-${\mathbbm B}$-isomorphic with $A$. This algebra ${\mathscr A}$
is unique up to isomorphism inside ${\mathbbm V}^{({\mathbbm B})}$.
\endproc

Observe that if ${\mathscr A}$ is an~$AW^\ast$-factor inside ${\mathbbm V}^{({\mathbbm B})}$
then the bounded descent $A$ of ${\mathscr A}$ is an~$AW^\ast$-algebra
whose Boolean algebra of central projections is isomorphic with~$\mathbbm B$.
Conversely, if $A$ is an~$AW^\ast$-algebra and ${\mathbbm B}\!:={\mathfrak P}_c(A)$
then there is an $AW^\ast$-factor ${\mathscr A}$ inside~${\mathbbm V}^{({\mathbbm B})}$ whose
bounded descent is isomorphic with $A$.

\subsection{11.3}
Take an $AW^\ast$-algebra $A$. Clearly, the
formula
$$
q\leq p\leftrightarrow q=qp=pq
\quad(q, p\in{\mathfrak P}(X))
$$
(sometimes reads as ``$p$ contains $q$'') specifies some order
$\leq$ on the set of projections ${\mathfrak P}(A)$. Moreover,
${\mathfrak P}(A)$ is a~complete lattice and ${\mathfrak P}_c(A)$ is a~complete
Boolean algebra.

The classification of $AW^\ast$-algebras into types is determined
from the structure of its lattice of projections \cite{DOP, Sak}. It
is important to emphasize that Boolean valued representation preserves
this classification. We recall  only the definition of type~I
$AW^\ast$-algebra. A~projection $\pi\in A$ is called
{\it abelian\/} if the algebra $\pi A\pi$ is commutative. An
algebra $A$ has {\it type\/}~I, if each nonzero projection in
$A$ contains a~nonzero abelian projection.

We call an $AW^*$-algebra {\it embeddable\/} if it is $*$-isomorphic
with the double commutant of some type~I $AW^*$-algebra. Each
embeddable $AW^*$-algebra admits a~Boolean valued representation,
becoming a~von Neumann algebra or factor.
A~$C^\ast$-algebra $A$ is called
${\mathbbm B}$-{\it embeddable} if there is
a~type~I $AW^\ast$-algebra $N$ and  a~$\ast$-monomorphism
$\imath: A\rightarrow N$ such that ${\mathbbm B} ={\mathfrak P}_c(N)$
and
$\imath (A) =\imath (A)^{\prime\prime}$, where
$\imath (A)^{\prime\prime}$ is the bicommutant of~$\imath (A)$
in~$N$. Note that in this event $A$ is an~$AW^\ast$-algebra and
${\mathbbm B}$ is a~regular subalgebra of~${\mathfrak P}_c(A)$.
In particular, $A$ is a~${\mathbbm B}$-cyclic algebra
(see~10.3).

Say that a~$C^\ast$-algebra~$A$  is {\it embeddable\/} if $A$
is
${\mathbbm B}$-embeddable for some
regular subalgebra
${\mathbbm B}\subset{\mathfrak P}_c(A)$. If ${\mathbbm B} ={\mathfrak P}_c(A)$
and $A$ is ${\mathbbm B}$-embeddable then $A$ is called
a~{\it centrally embeddable algebra}.

\subsection{11.4}
\proclaim{Theorem.}
Let ${\mathscr A}$ be a~$C^\ast\!$-algebra inside~${{\mathbbm V}}^{({\mathbbm B})}$
and let $A$ be the bounded  descent of~${\mathscr A}$. Then $A$ is
a~${\mathbbm B}$-embeddable $AW^\ast\!$-algebra if and only if
${\mathscr A}$ is a~von Neumann algebra inside~${{\mathbbm V}}^{({\mathbbm B})}$.
The algebra $A$ is centrally embeddable if and only if
${\mathscr A}$ is a~von Neumann factor inside~${{\mathbbm V}}^{({\mathbbm B})}$.
\endproc

Using this representation, we can obtain characterizations
of embeddable $AW^*$-algebras. In particular, an $AW^*$-algebra $A$
is embeddable if and only if the center-valued normal states of $A$
separate~$A$.

\subsection{11.5}
\proclaim{Theorem.}
For an~$AW^\ast$-algebra
$A$ the following are equivalent:

\subsubsec{(1)}~$A$ is embeddable;

\subsubsec{(2)}~$A$ is centrally embeddable;

\subsubsec{(3)}~$A$ has a~separating set of center-valued
normal states;

\subsubsec{(4)}~$A$ is a~${\mathfrak P}_c(A)$-predual space.
\endproc

\subsection{11.6}
Combining the results about the Boolean
valued representations of $AW^\ast$-algebras with
the analytical representations for dominated operators (see \cite{DOP}),
we come to some functional representations of~$AW^\ast$-algebras.

Suppose that $Q$ is an~extremally disconnected compact space, $H$ is
a~Hilbert space, and $B(H)$ is the space of bounded linear
endomorphisms of~$H$. Denote by ${\mathfrak C} (Q,B(H))$
the set of all operator-functions $u:\dom(u)\rightarrow B(H)$
 on the comeager sets $\dom(u)\subset Q$ and continuous
in the strong  operator topology. Introduce some equivalence on
${\mathfrak C} (Q,B(H))$ by  putting $ u\sim v$ if and only if $u$
and $v$ agree on $\dom(u)\cap\dom(v)$.

If $u\in {\mathfrak C}(Q,B(H))$ and $h\in H$ then the
vector-function  $uh:q\mapsto u(q)h \; (q\in\dom(u))$
is continuous thus determining  a~unique element
$\widetilde{uh}\in C_\infty (Q,H)$ from the condition $uh\in\widetilde{uh}$.
If $\tilde{u}$~is the coset of the operator-function
$u:\dom(u)\rightarrow B(H)$ then
$\tilde{u}h\!:=\widetilde{uh}\enskip (h\in H)$
by definition.

Denote by $SC_\infty (Q,B(H))$
the set of all cosets $\tilde{u}$ such that $u\in{\mathfrak C} (Q,B(H))$
and the set $\{\[\tilde{u}h\]:\,\|h\|\leq{\mathbb 1}\}$ is bounded
in $C_\infty (Q)$. Put
$$
\[\tilde{u}\]:= \sup\{\[\tilde{u}h\]:\|h\|\leq{\mathbb 1}\},
$$
where the supremum is taken in $C_\infty (Q)$.

We
naturally furnish $SC_\infty (Q,B(H))$  with the structure
of a~$\ast$-algebra  and unitary $C_\infty (Q)$-module.
We now introduce the following normed $\ast$-algebra
$$
\allowdisplaybreaks
\gathered
 SC_{\#} (Q,B(H))\!:=\{v\in SC_{\infty} (Q,B(H)):\[ v\]
\in C(Q)\},
\\
\| v\|=\|\[ v\]\|_ {\infty} \quad (v\in SC_{\#} (Q,B(H))).
\endgathered
$$

\subsection{11.7}
\proclaim{Theorem.}
To each type I $AW^\ast$-algebra $A$ there exists
a~family of nonempty extremally disconnected compact spaces
$(Q_\gamma)_{\gamma\in\Gamma}$
such that

\subsubsec{(1)}
$\Gamma$~is a~set
of cardinals and  $Q_\gamma$ is $\gamma$-stable
for every $\gamma\in\Gamma$;

\subsubsec{(2)}
there is a~$\ast$-${\mathbbm B}$-isomorphism:

$$
A\simeq \sum_{\gamma\in\Gamma}\!{}^{{}^\oplus}
SC_\# (Q_\gamma, B(l_2 (\gamma))).
$$

\noindent
This family is unique
up to congruence.
\endproc

A cardinal number $\gamma$ is  $Q$-{\it stable\/} if
$\gamma^{\scriptscriptstyle\wedge}$~is
a cardinal number inside ${\mathbbm V}^{({\mathbbm B})}$ and $Q$ is the Stone space
of~${\mathbbm B}$.

\subsection{11.8}~The study of $C^\ast$-algebras and von Neumann
algebras by  Boolean valued models was started by G.~Takeuti
with~\cite{Tak3, Tak4}. Theorems 11.2, 11.4, and 11.5
were obtained by M.~Ozawa \cite{Ozawa3, Ozawa4, Ozawa6}. Theorem~11.7 was established by
A.~G.~Kusraev.

Boolean valued analysis of $AW^\ast$-algebras
yields a~negative solution to
the I.~Kaplansky problem of unique decomposition of a~type~I
$AW^\ast$-algebra into the direct sum of homogeneous bands.
M.~Ozawa gave this solution in~\cite{Ozawa4, Ozawa5}.
The lack of uniqueness is tied with the effect of the cardinal
shift that may happens on ascending  into a~Boolean valued
model ${\mathbbm V}^{({\mathbbm B})}$. The cardinal shift is impossible
in the case when the Boolean algebra of central idempotents
${\mathbbm B}$ under study satisfies the countable chain condition,
and so the decomposition in question is unique. I.~Kaplansky
established uniqueness of the decomposition on assuming that $B$
satisfies the countable chain condition and conjectured that
uniqueness fails in general~\cite{Kap3}.

\ssection{12. $JB$-Algebras}

We also consider similar problems for the so-called $JB$-algebras
presenting some real nonassociative analogs of $C^*$-algebras.

\subsection{12.1}
Recall that a~$JB$-algebra is a~real Banach
space $A$ which is  a~unital Jordan algebra
satisfying  the conditions:

\subsubsec{(1)}~$\|xy\|\le\|x\|\cdot\|y\|\quad (x, y\in A)$;

\subsubsec{(2)}~$\|x^2\|=\|x\|^2 \quad (x\in A)$;

\subsubsec{(3)}~$\|x^2\|\le\|x^2+y^2\|\quad(x, y\in A)$.

The intersection of all maximal associative subalgebras of~$A$
is called the {\it center} of $A$ and denoted by ${\mathscr Z}(A)$. Evidently
${\mathscr Z}(A)$ is an~associative $JB$-algebra and every such algebra
is isometrically isomorphic to the real Banach algebra $C(K)$ of
continuous functions on a~compact space~$K$. If
${\mathscr Z}(A)={\mathbbm R}\cdot{{\mathbb 1}}$ then $A$ is said to be
a~$JB$-{\it factor}.

The idempotents of $JB$-algebras are also called {\it projections\/}.
The set ${\mathfrak P}(A)$ of projections  naturally underlies a~complete
lattice.  The set ${\mathfrak P}_c(A)$ of all
projections belonging to the center makes a~Boolean algebra.
Assume that ${\mathbbm B}$ is a~subalgebra of~${\mathfrak P}_c(A)$. Then we say
that $A$ is a~${\mathbbm B}$-$JB$-{\it algebra\/} if
to each partition of unity $(e_\xi)_{\xi \in \Xi}$ in ${\mathbbm B}$
and  each family $(x_\xi)_{\xi \in \Xi}$ in $A$ there exists a
unique ${\mathbbm B}$-{\it mixing\/} $x:=\mix_{\xi \in \Xi}\,(e_\xi x_\xi)$;
i.e., a~unique $x\in A$ such that
$e_\xi x_\xi=e_\xi x$ for all $\xi\in\Xi$. If
${\mathbbm B}\,({\mathbbm R})\!:={\mathscr R}\!\downarrow={\mathscr Z}(A)$ then a ${\mathbbm B}$-$JB$-algebra
is also called a~{\it centrally extended $JB$-algebra.}

Clearly, the unit ball of
a~${\mathbbm B}$-$JB$-algebra $A$ is closed under
${\mathbbm B}$-mixing and so $A$ is a~${\mathbbm B}$-cyclic
Banach space. Therefore, from 9.4 we can arrive to

\subsection{12.2}
\proclaim{Theorem.}
The bounded descent of~a $JB$-algebra inside
${\mathbbm V}^{({\mathbbm B})}$ is a~${\mathbbm B}$-$JB$-algebra. Conversely, for every
${\mathbbm B}$-$JB$-algebra $A$ there exists a~unique
$JB$-algebra ${\mathscr A}$ (up to isomorphism)  whose bounded descent is isometrically
${\mathbbm B}$-isomorphic to~$A$. Moreover, $[\![{\mathscr A}$ is
a~$JB$-factor $\!]\!]={\mathbb 1}$ if and only if \
${\mathbbm B}\,({\mathbbm R})={\mathscr Z}(A)$.
\endproc

\subsection{12.3}
Let $A$ be a~${\mathbbm B}$-$JB$-algebra and
$\Lambda:={\mathbbm B}\,({\mathbbm R})$. An operator
$\Phi\in A^{\scriptscriptstyle\#}$ is called a~$\Lambda$-{\it valued
state\/} if $\Phi\ge 0$ and $\Phi({\mathbb 1})={\mathbb 1}$. A state
$\Phi$ is {\it normal\/} if $\Phi(x)=\olim\Phi(x_\alpha)$ for every increasing net
$(x_{\alpha})$ in $A$ with $x:=\sup x_\alpha$.

If ${\mathscr A}$ is a~Boolean valued representation of~$A$
then the ascent $\phi:=\Phi{\uparrow}$ of $\Phi$ is a~normal state
on~${\mathscr A}$.  Conversely, if $[\![\phi$
is a  normal state on~${\mathscr A}]\!]={\mathbb 1}$ then the restriction
to $A$ of~$\phi{\downarrow}$ is  a~$\Lambda$-valued
normal state. We now give a~characterization  of the
${\mathbbm B}$-$JB$-algebras that are ${\mathbbm B}$-dual spaces.

\subsection{12.4}
\proclaim{Theorem.}
For a~${\mathbbm B}$-$JB$-algebra $A$ the following are equivalent:

\subsubsec{(1)}~$A$ is a~${\mathbbm B}$-dual space;

\subsubsec{(2)}~$A$ is monotone complete and admits a~separating
set of $\Lambda$-valued normal states.

If one of these conditions holds then the part of
$A{\scriptscriptstyle\#}$ consisting of $o$-continuous
operators serves as a~${\mathbbm B}$-predual of~$A$.
\endproc

This  is just a~Boolean valued interpretation of the following
theorem by F.~W. Shultz \cite{Shultz}: {\sl a~$JB$-algebra $A$ is a~dual
Banach space if and only if $A$ is monotone complete and has a
separating family of normal states}.

\subsection{12.5}
An algebra $A$ satisfying one of the equivalent
conditions 12.4\,(1,2), is called a~${\mathbbm B}$-$JBW$-{\it algebra}. If,
moreover, ${\mathbbm B}$ coincides with the set of all central
projections then $A$ is  a~${\mathbbm B}$-$JBW$-{\it factor}.
It follows from  Theorem 12.4 that $A$ is
a~${\mathbbm B}$-$JBW$-algebra $({\mathbbm B}$-$JBW$-factor) if and only if
its Boolean valued representation ${\mathscr A}\in {\mathbbm V}^{({\mathbbm B})}$
is a~$JBW$-algebra ($JBW$-factor).

Consider one example. Let $X$ be a~Kaplansky--Hilbert module over
the algebra $\overline\Lambda:={\mathbbm B}\,({\mathbbm C}):={\mathscr C}\!\downarrow$. Then $X$ is
a~${\mathbbm B}$-cyclic Banach space and ${\mathscr L}_{\mathbbm B}(X)$ is
a~type~I $AW^*$-algebra.

Given $x,y\in X$, define the  seminorm
$$
p_{x,y}(a):=\|\langle ax,y\rangle\|_\infty\quad(a\in{\mathscr L}_{\mathbbm B}(X)),
$$
where $\langle\cdot,\cdot\rangle$ is the inner product on~$X$
with values in~$\overline\Lambda$. Denote by $\sigma_\infty$ the
topology on~${\mathscr L}_{\mathbbm B}(X)$ that is generated by the system
of seminorms $p_{x,y}$. Then a~$\sigma_\infty$-closed
${\mathbbm B}$-$JB$-algebra of selfadjoint operators presents an
example of a~${\mathbbm B}$-$JBW$-algebra.

\subsection{12.6}
Let $A$ be an associative algebra over a~field  of characteristic~${\neq 2}$.
Define the new multiplication $a\circ b:=1/2(ab+ba)$ on the vector space
of~$A$. Denote the resulting algebra by~$A^J$. This
$A^J$ is a~Jordan algebra. If the subspace $A_\circ$ of~$A$ is closed
under $a\circ b$ then $A_\circ$  is a~subalgebra of
$A^J$ and so $A_\circ$  is Jordan. Such a~Jordan algebra $A_\circ $
is called {\it special.} The nonspecial Jordan algebras
are referred to as {\it exceptional}.

Let ${{\mathbbm O}}$
be the {\it Cayley\/} or
{\it  octonian algebra}.
Let $M_n({\mathbbm O})$ be the algebra of
$n\times n$-matrices with entries in~${{\mathbbm O}}$. The involution~$*$
on~$M_n({\mathbbm O})$ is as usual the transposition of a~matrix
followed by conjugation of every entry.
The set $M_n({\mathbbm O})_{\sa}\!:=\{x\in M_n({\mathbbm O}):\,x^*=x\}$
of hermitian matrices is closed in~$M_n({\mathbbm O})$
under the Jordan multiplication
$x\circ y=1/2(xy+yx)$. The real vector space
$M_n({\mathbbm O})_{\sa}$  is a~Jordan algebra
under  $\circ$ only for  $n\leq 3$. In case  $n=1,2$ we arrive at
special Jordan algebras. The Jordan algebra
$M_3({\mathbbm O})_{\sa}$ is special and denoted by
$M_3^8$.

\subsection{12.7}
\proclaim{Theorem.}
A special ${\mathbbm B}$-$JB$-algebra $A$ is a~${\mathbbm B}$-$JBW$-algebra
if and only if $A$ is isomorphic to a~$\sigma_\infty$-closed
${\mathbbm B}$-$JB$-subalgebra of ${\mathscr L}_{\mathbbm B}(X)_{sa}$ for some
Kaplansky--Hilbert module $X$.
\endproc

\subsection{12.8}
\proclaim{Theorem.}
Each ${\mathbbm B}$-$JBW$-factor $A$ admits a~unique decomposition
$A=eA\oplus e^*A$ with a~central projection $e\in{\mathbbm B}$,
$e^*:={\mathbb 1}-e$, such that the algebra $eA$ is special and the
algebra $e^*A$ is purely exceptional. Moreover, $eA$ is
${\mathbbm B}$-isomorphic to a~$\sigma_\infty$-closed subalgebra of
selfadjoint endomorphisms of some $AW^*$-module and $e^*A$ is
isomorphic to $C(Q,M_3^8)$, where $Q$ is the Stone compact space
of the Boolean algebra $e^*{\mathbbm B}:=[0,e^*]$.
\endproc

\subsection{12.9}
The $JB$-algebras are   nonassociative real
analogs of $C^*$-algebras and von Neumann operator algebras.
The theory of these algebras stems from the article of
P.~Jordan, J.~von Neumann, and E.~Wigner \cite{JNW} and  exists as a~branch of functional
analysis  since the mid 1960s, when
D.~M.~Topping \cite{Top} and E.~St\o rmer \cite{Sto}
have started the study of the nonassociative real analogs
of von Neumann algebras, the $JW$-algebras presenting weakly closed
Jordan algebras of bounded  selfadjoint  operators in a~Hilbert space.
The steps of development are reflected in \cite{Aju, ASS, HOS}.
The Boolean valued approach to $JB$-algebras is outlined by
A.~G.~Kusraev.
More details and references are collected in~\cite{IBA}.

%\def\bibitem#1#2{\item[#2]}
%\def\bibitem#1#2#3{\bibitem[#3]}
%\markboth{\sc References}{\sc References}

%\ssection{}{References}
%\begin{center}{\bf\scshape References}\end{center}

\end{document}